\documentclass[reqno]{amsart}
\usepackage[margin=26truemm,tmargin=26truemm,lmargin=28truemm,rmargin=28truemm]{geometry}

\usepackage{amssymb}
\usepackage{amsthm}
\usepackage{mathtools}
\mathtoolsset{showonlyrefs, showmanualtags}

\usepackage{aligned-overset}

\usepackage{bm}
\usepackage{mleftright}
\usepackage{enumitem}

\numberwithin{equation}{section}

\theoremstyle{plain}
\newtheorem{theorem}{Theorem}[section]
\newtheorem{corollary}[theorem]{Corollary}
\newtheorem{proposition}[theorem]{Proposition}

\theoremstyle{definition}
\newtheorem{definition}{Definition}

\theoremstyle{remark}
\newtheorem{remark}{Remark}
\newtheorem{example}[theorem]{Example}

\renewcommand{\S}{\mathbb{S}}

\newcommand{\Laplacian}{\Delta}

\newcommand{\R}{\mathbb{R}}
\newcommand{\N}{\mathbb{N}}
\newcommand{\C}{\mathbb{C}}

\newcommand*{\variabledot}{\makebox[1ex]{\textbf{$\cdot$}}}

\newcommand{\Nullset}{\mathcal{N}}

\newcommand{\abs}[1]{\lvert#1\rvert} 
\newcommand{\Abs}[1]{\mleft \lvert #1 \mright \rvert}
\newcommand{\norm}[1]{\lVert#1\rVert} 
 
\newcommand{\set}[2]{\{ \, #1 : #2 \, \} }

\newcommand{\Diracopm}[1]{ H_{#1} }
\newcommand{\Hankelnu}[1]{\mathcal{H}_{#1}}

\newcommand{\Tnu}[1]{T_{#1}}
\newcommand{\DiracTnum}[2]{\widetilde{T}_{#1, #2}}

\newcommand{\Umunu}[2]{U_{#1, #2}}
\newcommand{\Kernelmunu}[2]{\mathcal{K}_{#1, #2}}

\newcommand{\const}{\bm{C}}

\newcommand{\Sconstd}[1]{\const_{\textup{S}}^{(#1)}}
\newcommand{\Sconstwpd}[3]{\const_{\textup{S}}^{(#3)}(#1, #2)}

\newcommand{\Diracconstwpmd}[4]{\const_{\textup{D}, #3}^{(#4)}(#1, #2)}

\newcommand{\BesselI}[1]{I_{#1}}
\newcommand{\BesselJ}[1]{J_{#1}}
\newcommand{\BesselK}[1]{K_{#1}}

\newcommand{\Cmunu}[2]{\mathrm{C}^{(#2)}_{#1}}
\newcommand{\Dmunu}[2]{\mathrm{D}^{(#2)}_{#1}}

\newcommand{\LegendreP}[2]{\mathrm{P}_{#1}^{#2}}
\newcommand{\LegendreQ}[2]{\mathrm{Q}_{#1}^{#2}} 

\newcommand{\hypergeometric}[4]{ \,{}_2 F_{1}\mleft( \begin{matrix}
		#1, #2 \\
		#3
	\end{matrix} ; #4 \mright)  }
\newcommand{\Pochan}[2]{(#1)_#2}

\newcommand{\gmns}{g_{\mu, \nu, s}}
\newcommand{\Lpnunu}[3]{L^{#1}_{#2, #3}}

\newcommand{\intervaloo}[2]{(#1, #2)}
\newcommand{\intervalco}[2]{[#1, #2)}
\newcommand{\intervaloc}[2]{(#1, #2]}

\usepackage[numbers,sort]{natbib}
\providecommand*{\doi}[1]{doi:\href{https://doi.org/#1}{#1}}
\renewcommand*{\MR}[1]{\href{https://mathscinet.ams.org/mathscinet-getitem?mr=#1}{MR#1}}

\usepackage{hyperref}
\hypersetup{
	setpagesize=false,
	bookmarksnumbered=true,
	bookmarksopen=true,
	colorlinks=true,
	linkcolor=blue,
	citecolor=blue,
}

\usepackage[T1]{fontenc}
\usepackage{lmodern}

\title[Integral transforms involving squares of the Bessel functions]{Identities and inequalities for integral transforms involving squares of the Bessel functions}
\date{}

\author{Soichiro Suzuki}
\address[Soichiro Suzuki]{Department of Mathematics, Chuo University, 1-13-27, Kasuga, Bunkyo-ku, Tokyo, 112-8551, Japan}
\email{soichiro.suzuki.m18020a@gmail.com}
\thanks{The author was supported by Japan Society for the Promotion of Science (JSPS) KAKENHI Grant Number JP23KJ1939.}


\begin{document}
	\begin{abstract}
		We consider an integral transform given by $T_{\nu} f(s) := \pi \int_0^\infty rs J_{\nu}(r s)^2 f(r) \, dr$, where $J_{\nu}$ denotes the Bessel function of the first kind of order $\nu$.
		As shown by \citeauthor{Wal2002} (\citeyear{Wal2002}), this transform plays an essential role in the study of optimal constants of smoothing estimates for the free Schr\"{o}dinger equation on $\mathbb{R}^d$. 
		On the other hand, \citeauthor{BSS2015} (\citeyear{BSS2015})  studied these optimal constants using a different method, and obtained a certain alternative expression for $T_{\nu} f$ involving the $d$-dimensional Fourier transform of $x \mapsto f(\lvert x \rvert)$ when $\nu = k + d/2 - 1$ for $k \in \mathbb{N}$. 
		In this paper, we extend their identity to non-integer indices and derive several inequalities from it. In particular, we give a dimension-comparison result for the smoothing estimates on $\mathbb{R}^d$ and $\mathbb{R}^{d+1}$.
	\end{abstract}
	\maketitle
	\section{Introduction}
	We consider an integral transform $\Tnu{\nu}$ given by
	\begin{equation}
		\Tnu{\nu}f(s) \coloneqq \pi \int_0^\infty rs \BesselJ{\nu}(r s)^2 f(r) \, dr ,
		\label{eq:Tnu} \noeqref{eq:Tnu}
	\end{equation}
	where $\BesselJ{\nu}$ denotes the Bessel function of the first kind of order $\nu \in \intervalco{-1/2}{\infty}$.
	It is easy to see that the right-hand side exists as a Lebesgue integral for every $s \in \intervaloo{0}{\infty}$ when 
	\begin{equation}
		\int_0^1 \abs{f(r)} r^{2\nu+1}\, dr + \int_1^\infty \abs{f(r)} \, dr < \infty
	\end{equation}
	holds by a standard estimate \cite[\href{https://dlmf.nist.gov/10.7.E3}{10.7.3}, \href{https://dlmf.nist.gov/10.7.E8}{10.7.8}]{DLMF}
	\begin{equation} \label{eq:BesselJ upper bound}
	\BesselJ{\nu}(r) = \begin{cases*}
	O(r^{\nu}) , & as $r \downarrow 0$, \\
	O(r^{-1/2}) , & as $r \to \infty$,
	\end{cases*}
	\end{equation}
	and in this case $\Tnu{\nu} f$ is continuous on $\intervaloo{0}{\infty}$. In what follows, 
	we write $\Lpnunu{1}{\nu_1}{\nu_2}$ for the space of all measurable functions $f \colon \intervaloo{0}{\infty} \to \C$ such that the following norm is finite:
	\begin{equation} \label{eq:Lpnunu} \noeqref{eq:Lpnunu}
		\norm{f}_{\Lpnunu{1}{\nu_1}{\nu_2}} \coloneqq 
		\int_0^1 \abs{f(r)} r^{\nu_1}\, dr + \int_1^\infty \abs{f(r)} r^{\nu_2} \, dr .
	\end{equation}
	The aims of this paper are to establish a certain identity for $\Tnu{\nu}$ (Theorem \ref{main thm:identity and inequalities}) and to derive several inequalities from it (Theorems \ref{main thm:d+1 vs d} and \ref{main thm:Tnuf for completely monotone f}).
	
	As shown by \citet{Wal2002}, the transform $\Tnu{\nu}$ plays an essential role in the study of the following inequality, known as the Kato--Yajima smoothing estimate for the free Schr\"{o}dinger equation on $\R^d$:
	\begin{equation} \label{eq:smoothing Schrodinger} 
		\int_{(x, t) \in \R^d \times \R}  w(\abs{x}) \abs{\psi(\abs{\nabla}) e^{it\Laplacian} u_0(x)}^2 \, dx \, dt \leq C\norm{u_0}^2_{L^2(\R^d)} ,
	\end{equation}
	or equivalently, 
\begin{equation} 
	\frac{1}{(2\pi)^d}\int_{t \in \R}\int_{x \in \R^d} w(\abs{x}) \Abs{ \int_{\xi \in \R^d} e^{ix \cdot \xi} \psi(\abs{\xi}) e^{ - i t \abs{\xi}^2 }  \widehat{u_0}(\xi) \, d\xi }^2 \, dx \, dt \leq  C \norm{\widehat{u_0}}^2_{L^2(\R^d)} .
\end{equation}
	Here, $\widehat{u_0}$ denotes the Fourier transform of the initial data $u_0 \in L^2(\R^d)$;
	\begin{equation}
	\widehat{u_0}(\xi) \coloneqq \frac{1}{(2\pi)^{d/2}} \int_{x \in \R^d} u_0(x) \exp(- i x \cdot \xi) \, dx, 
	\end{equation}
	and $w, \psi \colon \intervaloo{0}{\infty} \to \intervalco{0}{\infty}$ are some given non-negative functions, which will be referred to as a spatial weight and smoothing function, respectively. For example, it is classically known that the smoothing estimate \eqref{eq:smoothing Schrodinger} holds in the following cases:
	\begin{alignat}{8}
		&d \geq 3,  \quad && \quad && ( w(r), \psi(s) ) = ( &&(1+r^2)^{-1}, \quad&&(1+s^2)^{1/4} &&) ,
		\label{eq:type A} \tag{A}\\
		%	&d = 2,  \quad && a > 2 , \quad && ( w(r), \psi(s) ) = ( &&(1+r^2)^{-a/2}, \quad&&(1+s^2)^{1/4} &&) ,
		%	\label{eq:type A 2D} \tag{A2}\\
		&d \geq 2,  \quad && 1 < p < d, \quad&& ( w(r), \psi(s) ) = ( &&r^{-p}, \quad&&s^{(2-p)/2} &&) , 
		\label{eq:type B} \tag{B}\\
		&d \geq 2,  \quad&& p > 1 , \quad&& ( w(r), \psi(s) ) = ( &&(1+r^2)^{-p/2}, \quad&&s^{1/2} &&) .
		\label{eq:type C} \tag{C}
	\end{alignat}
	See \citet[Theorem 2]{KY1989} for \eqref{eq:type A}; 
	%\citet[Proposition 2]{BK1992} for \eqref{eq:type A 2D}, 
	\citet[Theorem 1, Remarks (a)]{KY1989}, \citet[Theorem 1.1]{Sug1998}, \citet[Theorem 3]{Wat1991} for \eqref{eq:type B}; \citet[Theorem 1.(b)]{BK1992}, \citet[Theorem 1.1]{Chi2002} for \eqref{eq:type C}.
	Furthermore, the exponents in these cases are sharp; 
	see \citet[Theorem 2.1.(b), Theorem 2.2.(b)]{Wal1999} for \eqref{eq:type A},
	% and \eqref{eq:type A 2D} 
	\citet[Theorem 2]{Vil2001} for \eqref{eq:type B}, and \citet[Theorem 2.14.(b)]{Wal2000} for \eqref{eq:type C}.  
	In what follows, we always assume that $w$ is measurable and $\psi$ is continuous.
	For each $(w, \psi)$ and dimension $d \in \N_{\geq 1}$, let $\Sconstwpd{w}{\psi}{d}$ be the optimal constant of the inequality \eqref{eq:smoothing Schrodinger}.
	\citet{Wal2002} proved the following formula for this optimal constant.
	\begin{theorem}[{\cite[Theorem 4.1]{Wal2002}}] \label{thm:Walther}
		Let $d \geq 2$. Then we have
		\begin{equation} \label{eq:Walther}
			\Sconstwpd{w}{\psi}{d} = \sup_{k \in \N} \sup_{s > 0} s^{-1} \psi(s)^2 \Tnu{k + d/2 - 1}w(s) .
		\end{equation}
	\end{theorem}
	Note that Theorem \ref{thm:Walther} immediately implies
	\begin{equation} \label{eq:d+2 vs d} \noeqref{eq:d+2 vs d}
		\Sconstwpd{w}{\psi}{d+2} \leq \Sconstwpd{w}{\psi}{d}
	\end{equation}
	for every pair $(w, \psi)$ whenever $d \geq 2$.
	Inspired by this observation, \citet[Theorem 3.3]{Suz2026_Bessel} showed that 
	\begin{equation} \label{eq:d+1 vs d with 2} \noeqref{eq:d+1 vs d with 2}
	\Sconstwpd{w}{\psi}{d+1} \leq 2 \Sconstwpd{w}{\psi}{d}
	\end{equation}
	holds for every $(w, \psi)$ and $d \geq 2$.
	On the other hand, it is known that the stronger estimate
	\begin{equation} \label{eq:d+1 vs d}
		\Sconstwpd{w}{\psi}{d+1} \leq \Sconstwpd{w}{\psi}{d}
	\end{equation}
	holds in the cases \eqref{eq:type A}, \eqref{eq:type B}, \eqref{eq:type C}.
	For instance, in the case \eqref{eq:type A}, we have:
	\begin{theorem}[{\cite[Theorem 1.7]{BS2017}}] \label{thm:type A} 
		In the case \eqref{eq:type A}, we have
		\begin{equation}
			\Sconstd{d} = \begin{cases}
				\pi , & d = 3 , \\
				\pi \sup_{s>0} (1+s^2)^{1/2} \BesselI{1}(s) \BesselK{1}(s) , & d = 4 , \\
				\pi / 2 , & d \geq 5 ,
			\end{cases} 
		\end{equation}
		where $\BesselI{\nu}$ and $\BesselK{\nu}$ are the modified Bessel functions of the first and second kinds of order $\nu$, respectively. 
		Here we note that 
		\begin{equation}
			\sup_{s>0} (1+s^2)^{1/2} \BesselI{1}(s) \BesselK{1}(s)
			= 
			0.50239\ldots .
		\end{equation}
	\end{theorem}
	See \cite[Corollary 4]{Wat1991}, \cite[Eq.~(3)]{Sim1992}, \cite[Theorem 1.6]{BS2017} for the case \eqref{eq:type B}; and see \cite[Eq.~(2)]{Sim1992}, \cite[Theorem 1.4, Corollary 1.5]{BSS2015} for the case \eqref{eq:type C}.
	Our first result is a sufficient condition for the inequality \eqref{eq:d+1 vs d}.
	\begin{theorem} \label{main thm:d+1 vs d}
		Let $d \geq 2$.
		Assume that both of the Fourier transforms of a spatial weight $w \in \Lpnunu{1}{d-1}{0}$ as radial functions on $\R^d$ and $\R^{d+1}$ are non-negative.
		Then we have
		\begin{equation} 
			\Sconstwpd{w}{\psi}{d+1} \leq \Sconstwpd{w}{\psi}{d}
		\end{equation}
		for every smoothing function $\psi \colon \intervaloo{0}{\infty} \to \intervalco{0}{\infty}$.
		Here, non-negativity of the Fourier transform on $\R^d$ (as well as on $\R^{d+1}$) should be understood in the sense that
		\begin{equation}
	\int_{x \in \R^d} w(\abs{x}) \exp(- \varepsilon \abs{x}^2) \exp(- i x \cdot \xi) \, dx \geq 0
\end{equation}
holds for every $\varepsilon \in \intervaloo{0}{\infty}$ and $\xi \in \R^d$. 
	\end{theorem}
	Note that the assumption of Theorem \ref{main thm:d+1 vs d} is satisfied in the cases \eqref{eq:type A}, \eqref{eq:type B}, \eqref{eq:type C}, since the Riesz and Bessel kernels are non-negative.
	
	Our proof of Theorem \ref{main thm:d+1 vs d} is inspired by \citet*{BSS2015}. 
	In order to explain our idea briefly, we introduce the Hankel transform.
	\begin{definition}[Hankel transform]
		Let $\nu \in [-1/2, \infty)$ and $f \colon \intervaloo{0}{\infty} \to \C$. The Hankel transform of $f$ of order $\nu$ is formally defined by
		\begin{equation}
			\Hankelnu{\nu}f(\rho) \coloneqq \rho^{-\nu} \int_0^\infty r^{\nu + 1} \BesselJ{\nu}(r \rho) f(r) \, dr .  \label{eq:Hankelnu}
		\end{equation}
		For example, the right-hand side exists as a Lebesgue integral for $f \in \Lpnunu{1}{2\nu+1}{\nu+1/2}$ due to \eqref{eq:BesselJ upper bound}, and in that case it is continuous on $\intervaloo{0}{\infty}$.
	\end{definition}
	We remark that when $d \in \N_{\geq 1}$, the Hankel transform $\Hankelnu{d/2 - 1} f$ coincides with the Fourier transform of $f$ as a radial function on $\R^d$ in the sense that
	\begin{equation}
		\Hankelnu{d/2 - 1} f(\abs{\xi}) = \frac{1}{(2\pi)^{d/2}} \int_{x \in \R^d} f(\abs{x}) e^{- i x \cdot \xi} \, dx 
	\end{equation}	
	holds whenever $f \in \Lpnunu{1}{d-1}{d-1}$ ($\iff f(\abs{\variabledot}) \in L^1(\R^d)$).
	
	Now we define another integral transform $\Umunu{k}{d/2 - 1}$ for each $k \in \N$ and $d \in \N_{\geq 2}$ as follows:
	\begin{gather}
	\Umunu{k}{d/2 - 1}f(s) \coloneqq \int_{r \in \intervaloo{0}{\infty}} \Kernelmunu{k}{d/2 - 1}(r, s) f(r) \, dr , \\
	\Kernelmunu{k}{d/2 - 1}(r, s) \coloneqq 
	\dfrac{\pi^{1/2} r^{d-2}}{2^{d/2 - 1} \Gamma( (d-1)/2 ) }
	\times
	\begin{cases}
		(1 - r^2/(4s^2))^{(d-3)/2} \Cmunu{k}{d/2 - 1}(1 - r^2/(2s^2))  , & 0 < r < 2s , \\
		0 , & r > 2s > 0 ,
	\end{cases}
\end{gather}
	where $\Cmunu{k}{d/2 - 1}$ denotes the Gegenbauer function of the first kind of degree $k$ and order $d/2 - 1$ (see Definition \ref{def:Gegenbauer}).
	\citet[Eq.~(1.11)]{BSS2015} showed that
	\begin{equation} \label{eq:BSS identity integer}
		\Tnu{k + d/2 - 1}f
		= \Umunu{k}{d/2 - 1} \Hankelnu{d/2 - 1} f
	\end{equation}
	holds whenever $k \in \N$, $d \in \N_{\geq 2}$, and $f \colon \intervaloo{0}{\infty} \to \C$ is sufficiently nice\footnote{The assumption given in \citet{BSS2015} is complicated and also slightly ambiguous (see \cite[page 1769]{BSS2015}). As we will see in Theorem \ref{main thm:identity and inequalities}, the identity \eqref{eq:BSS identity integer} holds for $f \in \Lpnunu{1}{d-1}{(d-1)/2}$.}. 
As presented in \citet[page 1771, Remarks~(1)]{BSS2015}, this identity is quite useful to handle the supremum with respect to $k \in \N$ in \citeauthor{Wal2002}'s formula \eqref{eq:Walther}.
In fact, we see that
		\begin{equation}
			\Tnu{k + d/2 - 1}f(s) \leq \Tnu{d/2 - 1}f(s)
		\end{equation}
		holds for every $k \in \N$, $d \in \N_{\geq 2}$, and $s \in \intervaloo{0}{\infty}$ if $\Hankelnu{d/2 - 1} f$ is non-negative, since
		\begin{equation}
	\abs{ \Cmunu{k}{d/2 - 1}(x) } \leq \Cmunu{k}{d/2 - 1}(1) = \Cmunu{0}{d/2 - 1}(x) = 1 
		\end{equation}
		holds for every $k \in \N$, $d \in \N_{\geq 2}$, and $x \in \intervaloc{-1}{1}$ (see \cite[(2.116)]{AH2012}).
This observation leads us to 
\begin{equation}
\Sconstwpd{w}{\psi}{d} = \sup_{k \in \N} \sup_{s > 0} s^{-1} \psi(s)^2 \Tnu{k + d/2 - 1}w(s) = \sup_{s > 0} s^{-1} \psi(s)^2 \Tnu{d/2 - 1}w(s)
\end{equation}
for $w$ whose Hankel transform is non-negative.
	Therefore, to prove Theorem \ref{main thm:d+1 vs d}, it suffices to show that 
	\begin{equation} \label{eq:d+1 vs d sufficient}
		\Tnu{d/2 - 1/2}w(s) \leq \Tnu{d/2 - 1}w(s)
	\end{equation}
	holds.
	We prove this inequality \eqref{eq:d+1 vs d sufficient} by extending the identity \eqref{eq:BSS identity integer} to non-integer indices.
	We define $\Umunu{\mu}{\nu}$ for non-integer indices as follows.
	\begin{definition}
		Let $\mu \in \intervalco{0}{\infty}$ and $\nu \in \intervaloo{-1/2}{\infty}$.
		We define an integral transform $\Umunu{\mu}{\nu}$ by
		\begin{gather}
	\Umunu{\mu}{\nu}f(s) \coloneqq \int_{r \in \intervaloo{0}{\infty}} \Kernelmunu{\mu}{\nu}(r, s) f(r) \, dr, 
	\label{eq:Umunu}\\
	\hspace{-7mm} \Kernelmunu{\mu}{\nu}(r, s) \coloneqq \dfrac{\pi^{1/2} r^{2\nu} }{2^{\nu} \Gamma(\nu+1/2) }
	\times
	\begin{cases}
		\phantom{- \sin{(\mu \pi)} } (1 - r^2/(4s^2))^{\nu - 1/2} \Cmunu{\mu}{\nu}(1 - r^2/(2s^2))  , & 0 < r < 2s , \\
		- \sin{(\mu \pi)} (r^2/(4s^2) - 1)^{\nu - 1/2} \Dmunu{\mu}{\nu}(r^2/(2s^2) - 1) , & r > 2s > 0 .
	\end{cases} \label{eq:kernel} 
\end{gather}
	Here, $\Cmunu{\mu}{\nu}$ and $\Dmunu{\mu}{\nu}$ denote the Gegenbauer functions of the first and second kinds of degree $\mu$ and order $\nu$, respectively (see Definition \ref{def:Gegenbauer}).
	\end{definition}
	Note that $\Kernelmunu{\mu}{\nu}$ does not vanish on the region $r > 2s > 0$ if $\mu \not \in \N$, which is slightly unexpected.
	Additionally, in order to state our results simply, we introduce the following convention.
	\begin{definition}[non-negativity and strict positivity of the Hankel transform]
		Let $\nu \in \intervalco{-1/2}{\infty}$ and $f \in \Lpnunu{1}{2\nu+1}{0}$. We define $f_\varepsilon \in \Lpnunu{1}{2\nu+1}{2\nu+1}$ by 
		\begin{equation}
			f_\varepsilon(r) \coloneqq f(r) \exp(- \varepsilon r^2)
		\end{equation}
		for each $\varepsilon \in \intervaloo{0}{\infty}$.
		\begin{itemize}
			\item We say that \emph{$\Hankelnu{\nu} f$ is non-negative} if $\Hankelnu{\nu} f_\varepsilon(r) \geq 0$ holds for every $\varepsilon \in \intervaloo{0}{\infty}$ and $r \in \intervaloo{0}{\infty}$.
			\item We say that \emph{$\Hankelnu{\nu} f$ is strictly positive} if it is non-negative and
			\begin{equation}
				\liminf_{\varepsilon \downarrow 0} \Hankelnu{\nu} f_\varepsilon(r) > 0
			\end{equation}
			holds for almost every $r \in \intervaloo{0}{\infty}$.
		\end{itemize}
	\end{definition}
	Now our second main result is stated as follows.
	\begin{theorem} \label{main thm:identity and inequalities}
		The following hold:
		\begin{enumerate}[label=\textup{(\Roman*)}]
			\item \label{item:identity 1}
			Let $\mu \in \intervalco{0}{\infty}$, $\nu \in \intervaloo{-1/2}{\infty}$, and $f \in \Lpnunu{1}{2\nu+1}{\nu+1/2}$.
			Then, for every $s \in \intervaloo{0}{\infty}$, 
			\begin{equation}
				\Tnu{\mu + \nu} f(s) = \Umunu{\mu}{\nu} \Hankelnu{\nu} f(s)
			\end{equation}
			holds.
			\item \label{item:inequality 2}
			Let $\nu \in \intervalco{1/2}{\infty}$ and $f \in \Lpnunu{1}{2\nu+1}{0}$. 
			Then
			\begin{equation}
				\intervaloo{0}{\infty} \ni s \longmapsto \Tnu{\nu} f(s)
			\end{equation}
			is non-decreasing if $\Hankelnu{\nu} f$ is non-negative, and strictly increasing if $\Hankelnu{\nu} f$ is strictly positive.
			\item \label{item:inequality 3}
			Let $\nu \in \intervaloo{-1/2}{\infty}$ and $f \in \Lpnunu{1}{2\nu+1}{0}$. Then 
			\begin{equation}
				\intervaloo{0}{\infty} \ni s \longmapsto \Tnu{\nu} f(s) + \Tnu{\nu+1} f(s)
			\end{equation}
			is non-decreasing if $\Hankelnu{\nu} f$ is non-negative, and strictly increasing if $\Hankelnu{\nu} f$ is strictly positive.
			\item \label{item:inequality 4}
			Let $\mu \in \N_{\geq 1}$, $\nu \in \intervalco{0}{\infty}$, and $f \in \Lpnunu{1}{2\nu+1}{0}$. 
			Then, for every $s \in \intervaloo{0}{\infty}$, 
			\begin{equation}
				\Tnu{\mu + \nu} f(s) \leq \Tnu{\nu} f(s)
			\end{equation}
			holds if $\Hankelnu{\nu} f$ is non-negative, and the inequality is strict if $\Hankelnu{\nu} f$ is strictly positive.
			\item \label{item:inequality 5}
			Let $\mu \in \intervaloc{0}{1}$ and $\nu \in \intervaloo{-1/2}{\infty}$ be such that $\mu + 2\nu \geq 0$, and let $f \in \Lpnunu{1}{2\nu+1}{0}$. 
			Then, for every $s \in \intervaloo{0}{\infty}$, 
			\begin{equation}
				\Tnu{\mu + \nu} f(s) \leq \Tnu{\nu} f(s)
			\end{equation}
			holds if $\Hankelnu{\nu} f$ is non-negative, and the inequality is strict if $\Hankelnu{\nu} f$ is strictly positive.
		\end{enumerate}
	\end{theorem}	
	Here we present a short proof of Theorem \ref{main thm:d+1 vs d} using Theorem \ref{main thm:identity and inequalities}.
	\begin{proof}[Proof of Theorem \ref{main thm:d+1 vs d}]
		Since $\Hankelnu{d/2 - 1} w$ is non-negative by the assumption, 
\ref{item:inequality 4} and \ref{item:inequality 5} of Theorem \ref{main thm:identity and inequalities} imply that
\begin{equation} \label{eq:Tnuw ineq}
	\Tnu{\mu + d/2 - 1} w(s) \leq \Tnu{d/2 - 1} w(s) 
\end{equation}
holds for every $\mu \in \intervaloo{0}{1} \cup \N$ and $s \in \intervaloo{0}{\infty}$.
Therefore, Theorem \ref{thm:Walther} gives us
\begin{equation}
	\Sconstwpd{w}{\psi}{d}
	\underset{\text{Theorem \ref{thm:Walther}}}= \sup_{k \in \N} \sup_{s > 0} s^{-1} \psi(s)^2 \Tnu{k + d/2 - 1}w(s) 
	\underset{\eqref{eq:Tnuw ineq}}{=} \sup_{s > 0} s^{-1} \psi(s)^2 \Tnu{d/2 - 1} w(s) .
\end{equation}
By the same argument using the non-negativity of $\Hankelnu{d/2 - 1/2} w$, we also get
\begin{equation}
\Sconstwpd{w}{\psi}{d+1} = \sup_{s > 0} s^{-1} \psi(s)^2 \Tnu{d/2 - 1/2} w(s) .
\end{equation}
Now we use \eqref{eq:Tnuw ineq} with $\mu = 1/2$ and conclude that  
\begin{equation}
\Sconstwpd{w}{\psi}{d+1} \leq \Sconstwpd{w}{\psi}{d} 
\end{equation}
holds, as desired.
	\end{proof}
	Theorem \ref{main thm:identity and inequalities} is also useful to obtain monotonicity results.
	In particular, \ref{item:inequality 5} can be used to prove monotonicity of $\nu \longmapsto \Tnu{\nu} f$ by applying it to a function $f$ such that $\Hankelnu{\nu} f$ is non-negative for every $\nu$ (e.g., the Gaussian).
	Such a function is characterized in terms of the \emph{complete monotonicity}.
	\begin{definition}[complete monotonicity]
		A function $f \colon \intervaloo{0}{\infty} \to \R$ is said to be \emph{completely monotone} if it is infinitely differentiable and
		satisfies
		\begin{equation}
			(-1)^n f^{(n)}(r) \geq 0
		\end{equation}
		for every $n \in \N$ and $r \in \intervaloo{0}{\infty}$, where $f^{(n)}$ denotes the $n$-th derivative of $f$.
	\end{definition}
	\begin{theorem}[{\cite{Ber1929, Wid1931, Sch1938}}] \label{thm:completely monotone functions}
		Let $f \colon \intervaloo{0}{\infty} \to \R$. We consider the following conditions:
		\begin{enumerate}[label=\textup{(CM\arabic*)}]
			\item \label{item:completely monotone}
			The function $r \longmapsto f(r^{1/2})$ is completely monotone.
			\item \label{item:Gaussian mixture}
			There exists a Borel measure $\lambda$ on $\intervalco{0}{\infty}$ such that
			\begin{equation}
				f(r) = \int_{a \in \intervalco{0}{\infty}} e^{-a^2 r^2 / 2} \, d\lambda(a)
			\end{equation}
			holds for every $r \in \intervaloo{0}{\infty}$.
			\item \label{item:Hankel transform of measure for every order}
			For each $\nu \in \intervalco{-1/2}{\infty}$, there exists a Borel measure $\lambda_\nu$ on $\intervalco{0}{\infty}$ such that
			\begin{equation}
				f(r) = \Hankelnu{\nu} \lambda_\nu(r) \coloneqq \int_0^\infty (ar)^{-\nu} \BesselJ{\nu}(a r) \, d\lambda_\nu(a)
			\end{equation}
			holds for every $r \in \intervaloo{0}{\infty}$.
		\end{enumerate}
		Then we have
		\begin{equation}
			\ref{item:completely monotone} \iff \ref{item:Gaussian mixture} .
		\end{equation}
		This is known as the \citeauthor{Ber1929}--\citeauthor{Wid1931} theorem (\cite[\textsection 14]{Ber1929}, \cite[Theorem 8]{Wid1931}).
		Furthermore, if $f$ is bounded, then
		\begin{equation}	
			\ref{item:Gaussian mixture} \iff \ref{item:Hankel transform of measure for every order}.
		\end{equation}
		This is known as the \citeauthor{Sch1938} theorem (\cite[Theorem 2]{Sch1938}).
	\end{theorem}
	Taking into account Theorem \ref{thm:completely monotone functions}, it is natural to consider $\Tnu{\nu} f$ for a function $f$ such that $r \mapsto f(r^{1/2})$ is completely monotone. 
	Our last main result is as follows.
	\begin{theorem} \label{main thm:Tnuf for completely monotone f}
		Let $\nu_0 \in \intervaloo{-1/2}{\infty}$, and let $f \in \Lpnunu{1}{2\nu_0+1}{0}$ be such that 
		\begin{equation}
			r \longmapsto f(r^{1/2})
		\end{equation}
		is completely monotone and not identically zero.
		Then the following hold:
		\begin{enumerate}[label=\textup{(\arabic*)}]
			\item \label{item:Tnuf for completely monotone 1}
			For every $\nu \in \intervalco{\nu_0}{\infty}$, 
			\begin{equation}
				s \longmapsto s^{-\nu-1/2} \Tnu{\nu} f(s^{1/2})
			\end{equation}
			is completely monotone on $\intervaloo{0}{\infty}$.
			\item \label{item:Tnuf for completely monotone 2}
			For each fixed $\nu \in \intervalco{ 1/2 }{\infty} \cap \intervalco{\nu_0}{\infty}$, 
			\begin{equation}
				s \longmapsto \Tnu{\nu}f(s)
			\end{equation}
			is strictly increasing on $\intervaloo{0}{\infty}$.
			\item \label{item:Tnuf for completely monotone 3}
			For each fixed $\nu \in \intervaloo{ -1/2 }{\infty} \cap \intervalco{\nu_0}{\infty}$, 
			\begin{equation}
				s \longmapsto \Tnu{\nu}f(s) + \Tnu{\nu+1}f(s)
			\end{equation}
			is strictly increasing on $\intervaloo{0}{\infty}$.	
			\item \label{item:Tnuf for completely monotone 4}
			For each fixed $s \in \intervaloo{0}{\infty}$, 
			\begin{equation}
				\nu \longmapsto \Tnu{\nu} f(s)
			\end{equation}
			is strictly decreasing on $\intervalco{0}{\infty} \cap \intervalco{\nu_0}{\infty}$.
			\item \label{item:Tnuf for completely monotone 5}
			For each fixed $\nu \in \intervaloo{-1/2}{0} \cap \intervalco{\nu_0}{\infty}$ and $s \in \intervaloo{0}{\infty}$, 
			\begin{equation}
				\Tnu{-\nu} f(s) < \Tnu{\nu} f(s)
			\end{equation}
			holds.
		\end{enumerate}
	\end{theorem}
	We note that \ref{item:Tnuf for completely monotone 1} can also be regarded as a sufficient condition for the analyticity of $\Tnu{\nu} f$, since completely monotone functions are analytic owing to \citeauthor{Ber1914}'s little theorem (\cite[\textsection 1]{Ber1929}). 
	In other words, we have:
	\begin{corollary} \label{cor:analytic}
		Let $\nu_0 \in \intervaloo{-1/2}{\infty}$, and let $f \in \Lpnunu{1}{2\nu_0+1}{0}$ be such that 
		\begin{equation}
			r \longmapsto f(r^{1/2})
		\end{equation}
		is completely monotone.
		Then $\Tnu{\nu} f$ is analytic on $\intervaloo{0}{\infty}$ for every $\nu \in \intervalco{\nu_0}{\infty}$.
	\end{corollary}
	As \citet{BS2017} pointed out, the analyticity of $\Tnu{\nu} f$ can be used to show the non-existence of extremisers for the smoothing estimate \eqref{eq:smoothing Schrodinger} (see \cite[Theorem 1.2]{BS2017} for details).
	They established the following sufficient condition for the analyticity, and showed that $\Tnu{k/2} f$ is analytic when $f(r) = (1+r^2)^{-p/2}$ for every $p \in \intervaloo{1}{\infty}$ and $k \in \N$. 
	\begin{theorem}[{\cite[Theorem 1.3]{BS2017}}] \label{thm:analytic by BS2017}
		Let $f \colon \intervaloo{0}{\infty} \to \C$ be infinitely differentiable. Suppose that there exists a constant $C > 0$ such that
		\begin{equation}
			\int_0^\infty r^n \abs{ f^{(n)}(r) } \, dr \leq C^{n+1} n!
		\end{equation}
		holds for each $n \in \N$. 
		In addition, we also assume that 
		\begin{equation}
			\int_0^\infty r^n \sup_{\theta \in I}{\abs{ f^{(n)}(\theta r) }} \, dr < \infty
		\end{equation}
		holds for every bounded closed interval $I \subset \intervaloo{0}{\infty}$ and $n \in \N$.
		Then $\Tnu{k/2} f$ is analytic on $\intervaloo{0}{\infty}$ for every $k \in \N$.
	\end{theorem}
	Certainly, the analyticity of $\Tnu{\nu} f$ for $f(r) = (1+r^2)^{-p/2}$ can be proved much more easily by using Corollary \ref{cor:analytic} instead of Theorem \ref{thm:analytic by BS2017}.
	In fact, it is clear that $f \in \Lpnunu{1}{0}{0} = L^1$ and that
	\begin{equation}
		r \longmapsto f(r^{1/2}) = (1+r)^{-p/2}
	\end{equation}
	is completely monotone.
	See \cite[Lemma 4.4]{BS2017} for the proof via Theorem \ref{thm:analytic by BS2017}.
	
	\subsection*{Organization of the paper}
	In Section \ref{section:Gegenbauer}, we introduce some basic facts on the Gegenbauer functions and the Hankel transform, which will be used in the following sections. 
	In Section \ref{section:proof of identity and inequalities}, we prove Theorem \ref{main thm:identity and inequalities}. 
	The identity \ref{item:identity 1} is proved in Subsection \ref{subsection:proof of the identity}, and the inequalities \ref{item:inequality 2}, \ref{item:inequality 3}, \ref{item:inequality 4}, \ref{item:inequality 5} are in Subsection \ref{subsection:proofs of the inequalities}.
	In Section \ref{section:completely monotone}, we prove Theorem \ref{main thm:Tnuf for completely monotone f}.
	In Section \ref{section:remarks}, we give some remarks.
	Analogues of Theorems \ref{main thm:identity and inequalities} and \ref{main thm:Tnuf for completely monotone f} for the case $\nu = -1/2$ are given in Subsection \ref{subsection:v=-1/2}, 
	and an analogue of Theorem \ref{main thm:d+1 vs d} for the Dirac equations is in Subsection \ref{subsection:Dirac}.
	In Subsection \ref{subsection:differentiation}, we show that \ref{item:inequality 3} can also be derived from \ref{item:inequality 5} rather than \ref{item:identity 1}.
	\section{Preliminaries} \label{section:Gegenbauer}
	\subsection{Gegenbauer functions}
	First, we give the definition of the Gegenbauer functions.
 		\begin{definition}[Gegenbauer functions] \label{def:Gegenbauer}
		Let $\mu \in [0, \infty)$ and $\nu \in \intervaloo{-1/2}{\infty}$. The Gegenbauer function of the first kind $\Cmunu{\mu}{\nu} \colon \intervaloc{-1}{1} \to \R$ and that of the second kind $\Dmunu{\mu}{\nu} \colon \intervaloo{1}{\infty} \to \R$ are defined by
		\begin{gather} 
			\Cmunu{\mu}{\nu}(x) \coloneqq \hypergeometric{-\mu}{\mu + 2\nu}{\nu+1/2}{(1-x)/2} , 
			\label{eq:normalized Genenbauer 1st} \\
						\Dmunu{\mu}{\nu}(x) \coloneqq \frac{ \Gamma(\mu+1) \Gamma(\nu+1/2) }{ \pi^{1/2} 2^{\mu} \Gamma(\mu+\nu+1) } x^{-(\mu + 2\nu)} \hypergeometric{\mu/2 + \nu}{(\mu + 2\nu + 1)/2}{\mu + \nu + 1}{1/x^2} ,
			\label{eq:normalized Genenbauer 2nd}
		\end{gather}
		respectively. 
		Here, $\,{}_2 F_{1}$ denotes the hypergeometric function, which is given by
		\begin{gather}
			\hypergeometric{a}{b}{c}{x} \coloneqq \sum_{n=0}^{\infty} \frac{\Pochan{a}{n} \Pochan{b}{n}}{\Pochan{c}{n}} \frac{x^n}{n!} , \\
			\Pochan{a}{n} \coloneqq \prod_{m=0}^{n-1} (a+m) = \frac{ \Gamma(a+n) }{\Gamma(a)} .
		\end{gather}
	\end{definition}
	We note that $\Cmunu{k}{\nu}$ becomes a polynomial of degree $k$ for each $k \in \N$. In that case, it is usually called the Gegenbauer polynomial. In particular, it coincides with the Chebyshev and Legendre polynomials when $\nu = 0, 1/2$, respectively.
	The following simple bounds are elementary and well known.
	\begin{proposition}
		\phantom{a}
	\begin{itemize}
		\item 
		Let $k \in \N$ and $\nu \in \intervalco{0}{\infty}$. Then we have
		\begin{equation} \label{eq:Cmunu upper bound integer}
			\abs{ \Cmunu{k}{\nu}(x) } \leq 1 
		\end{equation}
		for every $x \in \intervaloc{-1}{1}$.
		\item Let $\mu \in \intervaloc{0}{1}$ and $\nu \in \intervaloo{-1/2}{\infty}$ be such that $\mu + 2 \nu \geq 0$. Then we have
		\begin{equation} \label{eq:Cmunu upper bound}
		\Cmunu{\mu}{\nu}(x) \leq 1 
		\end{equation}
		for every $x \in \intervaloo{-1}{1}$. Moreover, this becomes equality if and only if $\mu + 2\nu = 0$.
		\item Let $\mu \in \intervalco{0}{\infty}$ and $\nu \in \intervaloo{-1/2}{\infty}$ be such that $\mu + 2 \nu \geq 0$. Then we have
		\begin{equation} \label{eq:Dmunu positive}
		\Dmunu{\mu}{\nu}(x) > 0
		\end{equation}
		for every $x \in \intervaloo{1}{\infty}$.
	\end{itemize}
	\end{proposition}
	For the first inequality \eqref{eq:Cmunu upper bound integer}, see \cite[(2.116)]{AH2012}, \cite[Theorem 7.32.1]{Sze1975},  \cite[\href{https://dlmf.nist.gov/18.14.E4}{18.14.4}]{DLMF}, for example. 
	To see \eqref{eq:Cmunu upper bound}, 
	notice that we have
	\begin{equation}
		\frac{\Pochan{-\mu}{n} \Pochan{\mu + 2\nu}{n}}{\Pochan{\nu+1/2}{n}} \leq 0
	\end{equation}
	for every $n \in \N_{\geq 1}$ under the assumption. Hence, 
	\begin{equation}
		\Cmunu{\mu}{\nu}(x) = \hypergeometric{-\mu}{\mu+2\nu}{\nu+1/2}{(1-x)/2} = \sum_{n = 0}^{\infty} \frac{\Pochan{-\mu}{n} \Pochan{\mu + 2\nu}{n}}{\Pochan{\nu+1/2}{n}} \frac{( (1-x)/2 )^n}{n!} \leq 1
	\end{equation}
	holds for every $x \in \intervaloo{-1}{1}$, and the equality holds if and only if $\mu + 2\nu = 0$.
	Similarly, \eqref{eq:Dmunu positive} follows from the fact that
	\begin{equation}
		\frac{\Pochan{\mu/2 + \nu}{n} \Pochan{(\mu + 2\nu + 1)/2}{n}}{\Pochan{\mu+\nu+1}{n}} \geq 0
	\end{equation}
	holds for every $n \in \N$ under the assumption.
	
	We also have the following asymptotic formulas for the Gegenbauer functions, which are special cases of more general results for the hypergeometric functions (see \cite[\href{https://dlmf.nist.gov/15.4.ii}{15.4.ii}]{DLMF}, for example).
		\begin{proposition} \label{prop:behavior of CD}
		Let $\mu \in \intervalco{0}{\infty}$ and $\nu \in \intervaloo{-1/2}{\infty}$. 
		Then we have 
		\begin{align}
			\Cmunu{\mu}{\nu}(x)
			&\begin{dcases}
				\to (-1)^{\mu} , & \mu \in \N, \\
				\to \frac{\cos{(\mu + \nu)\pi}}{ \cos{\nu \pi} } , & \nu \in \intervaloo{-1/2}{1/2} \text{ and } \mu \not \in \N , \\
				\sim - \frac{\sin{\mu \pi}}{\pi} \log{\mleft(\frac{2}{1+x}\mright)} , & \nu = 1/2 \text{ and } \mu \not \in \N , \\
				\sim \frac{ \Gamma(\nu+1/2) \Gamma(\nu - 1/2) }{ \Gamma(-\mu) \Gamma(\mu + 2\nu) } \mleft( \frac{2}{1+x} \mright)^{\nu-1/2} , & \nu \in (1/2,  \infty) \text{ and } \mu \not \in \N 
			\end{dcases}
			\shortintertext{as $x \downarrow -1$,} 
			\Cmunu{\mu}{\nu}(x)
			&\to 1
			\shortintertext{as $x \uparrow 1$,} 
			\Dmunu{\mu}{\nu}(x)
			& \begin{dcases}
				\to \dfrac{1}{\cos{\nu \pi}} 	, & \nu \in \intervaloo{-1/2}{1/2} , \\
				\sim \frac{1}{\pi} \log{ \mleft( \frac{2}{x-1} \mright) } , & \nu = 1/2 , \\
				\sim \frac{\Gamma(\mu+1) \Gamma(\nu+1/2) \Gamma(\nu - 1/2)  }{ \pi \Gamma(\mu + 2\nu) } \mleft( \frac{2}{x-1} \mright)^{\nu - 1/2} , & \nu \in \intervaloo{1/2}{\infty}
			\end{dcases}
			\shortintertext{as $x \downarrow 1$, and} 
			\Dmunu{\mu}{\nu}(x)
			&\sim \frac{ \Gamma(\mu+1) \Gamma(\nu+1/2) }{ \pi^{1/2} 2^{\mu} \Gamma(\mu+\nu+1) } x^{-(\mu + 2\nu)}
		\end{align}
		as $x \uparrow \infty$.
		Here, $A(x) \sim B(x)$ as $x \downarrow x_0$ means that $\lim_{x \downarrow x_0} A(x) / B(x) = 1$ holds (as well as for $x \uparrow x_0$). 
	\end{proposition}
	Finally, we remark the following connection formulas between the Gegenbauer and Legendre functions.
	\begin{definition}[Legendre functions]
	Let $\mu \in \intervaloo{-\infty}{1}$ and $\nu \in \intervalco{-\mu}{\infty}$. The Legendre function of the first kind $\LegendreP{\nu}{\mu} \colon \intervaloo{-1}{1} \to \R$ and that of the second kind $\LegendreQ{\nu}{\mu} \colon \intervaloo{1}{\infty} \to \R$ are defined by
	\begin{gather}
		\LegendreP{\nu}{\mu}(x) \coloneqq \frac{1}{\Gamma(1-\mu)} \mleft( \frac{1+x}{1-x} \mright)^{\mu/2} \hypergeometric{\nu+1}{-\nu}{1-\mu}{(1-x)/2} ,
		\label{eq:Legendre 1st} \\
		\LegendreQ{\nu}{\mu}(x) \coloneqq \frac{\Gamma(\mu+\nu+1)}{2^{\nu} \pi^{1/2} \Gamma(\nu+3/2)}  \frac{ (x^2 - 1)^{\mu / 2} }{ x^{\mu + \nu + 1} } \hypergeometric{(\mu+\nu+2)/2}{(\mu+\nu+1)/2}{\nu+3/2}{1/x^2} ,
		\label{eq:Legendre 2nd}
	\end{gather}
	respectively. 
\end{definition}
\begin{proposition} \label{prop:Gegenbauer and Legendre}
	Let $\mu \in \intervalco{0}{\infty}$ and $\nu \in \intervaloo{-1/2}{\infty}$. Then we have
		\begin{align}
	\LegendreP{\mu + \nu - 1/2}{-\nu+1/2}(x) 
	&= \frac{2^{-\nu + 1/2}}{\Gamma(\nu+1/2) } (1-x^2)^{(\nu - 1/2)/2} \Cmunu{\mu}{\nu}(x) 
	\label{eq:P sim C} 
	\shortintertext{for every $x \in \intervaloo{-1}{1}$, and}
	\LegendreQ{\mu + \nu - 1/2}{-\nu+1/2}(x) 
	&= \frac{ 2^{-\nu+1/2} }{ \Gamma(\nu+1/2) } (x^2 - 1)^{(\nu-1/2) / 2}  \Dmunu{\mu}{\nu}(x) 
	\label{eq:Q sim D}
\end{align}
for every $x \in \intervaloo{1}{\infty}$.
\end{proposition}
These connection formulas are special cases of Euler's transformation formula for the hypergeometric function (\cite[\href{https://dlmf.nist.gov/15.8.E1}{15.8.1}]{DLMF}, \cite[9.131.1]{GR2014}):
\begin{equation}
	\hypergeometric{a}{b}{c}{z} = (1-z)^{c - (a+b)} \hypergeometric{c-a}{c-b}{c}{z} .
\end{equation}
	\subsection{Parseval's identity for the Hankel transform}
	We present some sufficient conditions for Parseval's identity for the Hankel transform.
	One of the simplest cases is the following. 	
\begin{proposition} \label{prop:Parseval 1}
Let $\nu \in \intervalco{-1/2}{\infty}$,
let $f \in \Lpnunu{1}{2\nu+1}{\nu+1/2}$, and let $g \in \Lpnunu{1}{\nu+1/2}{2\nu+1}$. Then we have
	\begin{equation} \label{eq:Parseval duality form}
	\int_{r \in \intervaloo{0}{\infty}} f(r) \Hankelnu{\nu} g(r) r^{2\nu+1} \, dr = \int_{\rho \in \intervaloo{0}{\infty}} \Hankelnu{\nu} f(\rho) g(\rho) \rho^{2\nu+1} \, d\rho , 
\end{equation}
where both sides exist as Lebesgue integrals.
\end{proposition}
This can be shown easily via the standard estimate \eqref{eq:BesselJ upper bound} and Fubini's theorem, so that we omit details. 
In order to deal with Hankel transforms which exist as improper integrals but not as Lebesgue integrals, one can use the following result.
\begin{proposition}[{\citet[Theorem IV]{Mac1939}}] \label{prop:Parseval 2}
	Let $\nu \in \intervalco{-1/2}{\infty}$,
	let $f \in \Lpnunu{1}{\nu+1/2}{\nu+1/2}$, and let $g \in \Lpnunu{1}{\nu+1/2}{\nu-3/2}$ be such that $r^{\nu+1/2} g(r) \to 0$ as $r \to \infty$. 
	Assume that $\Hankelnu{\nu} g(\rho)$ exists as an improper integral for every but finitely many $\rho \in \intervaloo{0}{\infty}$, and that $\Hankelnu{\nu} g \in \Lpnunu{1}{\nu+1/2}{\nu+1/2}$.
	Then we have
	\begin{equation}
		\int_{r \in \intervaloo{0}{\infty}} f(r) g(r) r^{2\nu+1} \, dr = \int_{\rho \in \intervaloo{0}{\infty}} \Hankelnu{\nu} f(\rho) \Hankelnu{\nu} g(\rho) \rho^{2\nu+1} \, d\rho ,
	\end{equation}
	where both sides exist as Lebesgue integrals.
\end{proposition}
Observe that the assumption $f, \Hankelnu{\nu} g \in \Lpnunu{1}{\nu+1/2}{\nu+1/2}$ in Proposition \ref{prop:Parseval 2} can be replaced with $f \in \Lpnunu{1}{2\nu+1}{\nu+1/2}$ and $\Hankelnu{\nu} g \in \Lpnunu{1}{\nu+1/2}{2\nu+1}$ by using Proposition \ref{prop:Parseval 1} as follows. 
\begin{proposition} \label{prop:Parseval 3}
	Let $\nu \in \intervalco{-1/2}{\infty}$,
let $f \in \Lpnunu{1}{2\nu+1}{\nu+1/2}$, and let $g \in \Lpnunu{1}{\nu+1/2}{\nu-3/2}$ be such that $r^{\nu+1/2} g(r) \to 0$ as $r \to \infty$. 
Assume that $\Hankelnu{\nu} g(\rho)$ exists as an improper integral for every but finitely many $\rho \in \intervaloo{0}{\infty}$, and that $\Hankelnu{\nu} g \in \Lpnunu{1}{\nu+1/2}{2\nu+1}$.
Then we have
\begin{equation}
	\int_{r \in \intervaloo{0}{\infty}} f(r) g(r) r^{2\nu+1} \, dr = \int_{\rho \in \intervaloo{0}{\infty}} \Hankelnu{\nu} f(\rho) \Hankelnu{\nu} g(\rho) \rho^{2\nu+1} \, d\rho ,
\end{equation}
where both sides exist as Lebesgue integrals.
\end{proposition}
Proposition \ref{prop:Parseval 3} plays an essential role in the next section. 
\begin{proof}[Proof of Proposition \ref{prop:Parseval 3}] 
We write $G \coloneqq \Hankelnu{\nu} g$ for simplicity.
Then we have
	\begin{equation}
		\int_{r \in \intervaloo{0}{\infty}} f(r) \Hankelnu{\nu} G(r) r^{2\nu+1} \, dr = \int_{\rho \in \intervaloo{0}{\infty}} \Hankelnu{\nu} f(\rho) \Hankelnu{\nu} g(\rho) \rho^{2\nu+1} \, d\rho ,
	\end{equation}
	and both sides exist as Lebesgue integrals
	thanks to Proposition \ref{prop:Parseval 1}. 
Therefore, it suffices to show that
\begin{equation} \label{eq:inverse}
\Hankelnu{\nu} G(r) = g(r)
\end{equation}
holds for almost every $r \in \intervaloo{0}{\infty}$.
Now fix $\varphi \in C_c^\infty(\intervaloo{0}{\infty})$ arbitrarily.
Then Propositions \ref{prop:Parseval 1} and \ref{prop:Parseval 2} imply
\begin{align}
\int_{r \in \intervaloo{0}{\infty}} \varphi(r) \Hankelnu{\nu} G(r) r^{2\nu+1} \, dr &= \int_{\rho \in \intervaloo{0}{\infty}} \Hankelnu{\nu} \varphi(\rho) \Hankelnu{\nu} g(\rho) \rho^{2\nu+1} \, d\rho , \\
\int_{r \in \intervaloo{0}{\infty}} \varphi(r) g(r) r^{2\nu+1} \, dr &= \int_{\rho \in \intervaloo{0}{\infty}} \Hankelnu{\nu} \varphi(\rho) \Hankelnu{\nu} g(\rho) \rho^{2\nu+1} \, d\rho ,
\end{align}
respectively, so that 
\begin{equation}
\int_{r \in \intervaloo{0}{\infty}} \varphi(r) \Hankelnu{\nu} G(r) r^{2\nu+1} \, dr = \int_{r \in \intervaloo{0}{\infty}} \varphi(r) g(r) r^{2\nu+1} \, dr
\end{equation}
holds. Hence, we conclude that \eqref{eq:inverse} holds for almost every $r \in \intervaloo{0}{\infty}$ as desired, since $\varphi \in C_c^\infty(\intervaloo{0}{\infty})$ is arbitrary.
\end{proof}
	\section{\texorpdfstring{Proof of Theorem \ref{main thm:identity and inequalities}}{Proof of Theorem 1.4}} \label{section:proof of identity and inequalities}
	\subsection{\texorpdfstring{Proof of the identity \ref{item:identity 1}}{Proof of the identity (I)}} \label{subsection:proof of the identity}
	To begin with, we give a brief outline of our proof of \ref{item:identity 1}.
	We define $\gmns \colon \intervaloo{0}{\infty} \to \R$ by
	\begin{equation}
		\gmns(r) \coloneqq \pi s r^{-2 \nu} ( \BesselJ{\mu + \nu}(rs) )^2 ,
	\end{equation}
	so that
	\begin{equation} 
		\Tnu{\mu + \nu}f(s) = \int_{r \in \intervaloo{0}{\infty}} f(r) \gmns(r) r^{2\nu+1} \, dr .
	\end{equation}
	The first step is to find an explicit expression for the Hankel transform of $\gmns$. 
	More precisely, we will see that the kernel $\Kernelmunu{\mu}{\nu}$ defined in \eqref{eq:kernel}
	satisfies
	\begin{equation} \label{eq:kernel = Hankel transform}
		\rho^{2\nu+1} \Hankelnu{\nu} \gmns(\rho) = \Kernelmunu{\mu}{\nu}(\rho, s) .
	\end{equation}
	Now we recall Parseval's identity for the Hankel transform, which states that
	\begin{equation} \label{eq:Parseval}
		\int_{r \in \intervaloo{0}{\infty}} f(r) g(r) r^{2\nu+1} \, dr = \int_{\rho \in \intervaloo{0}{\infty}} \Hankelnu{\nu} f(\rho) \Hankelnu{\nu}g(\rho) \rho^{2\nu+1} \, d\rho 
	\end{equation}
	holds for suitable pairs of functions $(f, g)$.
	In the second step, we will show that Parseval's identity \eqref{eq:Parseval} is valid when $f \in \Lpnunu{1}{2\nu+1}{\nu + 1/2}$ and $g = \gmns$ by using Proposition \ref{prop:Parseval 3}.
	Then the desired identity \ref{item:identity 1} follows from \eqref{eq:kernel = Hankel transform} and \eqref{eq:Parseval}.
	
	Now we are going to prove the first step.
	\begin{proposition}\label{prop:explicit formula for kernel}
		We have
		\begin{equation} \label{eq:explicit formula for kernel}
			\Hankelnu{\nu} \gmns(\rho) 
			= \dfrac{\pi^{1/2}}{2^{\nu} \Gamma(\nu+1/2) \rho} \times 
			\begin{cases}
				\phantom{- \sin{(\mu \pi)} }  (1 - \rho^2/(4s^2))^{(2\nu - 1)/2} \Cmunu{\mu}{\nu}(1 - \rho^2/(2s^2))  , & 0 < \rho < 2s , \\
				-  \sin{(\mu \pi)} (\rho^2/(4s^2) - 1)^{(2\nu - 1)/2} \Dmunu{\mu}{\nu}(\rho^2/(2s^2) - 1) , & \rho > 2s > 0
			\end{cases}
		\end{equation}
		for every $\mu \in \intervalco{0}{\infty}$, $\nu \in \intervaloo{-1/2}{\infty}$, and $s \in \intervaloo{0}{\infty}$.
		Here, $\Hankelnu{\nu} \gmns(\rho)$ exists as a Lebesgue integral when $\nu \in \intervaloo{1/2}{\infty}$, but as an improper integral when $\nu \in \intervaloc{-1/2}{1/2}$.
	\end{proposition}
	Observe that this is just a special case of \citeauthor{Mac1909}'s formula (\cite{Mac1909}), which is usually stated using the Legendre functions instead of the Gegenbauer functions as follows.
	\begin{theorem}[{\cite{Mac1909}, \cite[\href{https://dlmf.nist.gov/10.22.E71}{10.22.71}, \href{https://dlmf.nist.gov/10.22.E72}{10.22.72}]{DLMF}, \cite[13.46.(4), (5)]{Wat1944}}] \label{thm:explicit formula Legendre ver}
		Let $a, b, c \in \intervaloo{0}{\infty}$, $\mu \in \intervalco{0}{\infty}$, and $\nu \in \intervaloo{-1/2}{\infty}$. Then we have
		\begin{align} 
			&\quad \int_{r \in \intervaloo{0}{\infty}} r^{-\nu+1} \BesselJ{\nu}(ar) \BesselJ{\mu + \nu}(br) \BesselJ{\mu + \nu}(cr) \, dr \\
			&= \mleft\{ 
			\begin{alignedat}{7}
				\phantom{- \frac{ 2 \sin{\mu \pi} }{ \pi }} \frac{ (bc)^{\nu - 1} (\sin{\phi})^{\nu - 1/2} }{ (2\pi)^{1/2} a^{\nu} } & \LegendreP{\mu + \nu - 1/2}{-\nu+1/2}(\cos{\phi}) , \quad \quad \abs{b-c} <{} &a < b + c , \\
				- \sin{(\mu \pi)} \frac{ (bc)^{\nu - 1} (\sinh{\chi})^{\nu - 1/2} }{ (2\pi)^{1/2} a^{\nu} } & \LegendreQ{\mu + \nu - 1/2}{-\nu+1/2}(\cosh{\chi}) , \quad &  a > b + c ,
			\end{alignedat} 
			\mright. \label{eq:explicit formula Legendre ver}
		\end{align}
		where $\phi \in (0, \pi)$ and $\chi \in \intervaloo{0}{\infty}$ are such that 
		\begin{equation}
			\cos{\phi} = \frac{ b^2 + c^2 - a^2 }{ 2bc } , \quad \cosh{\chi} = \frac{a^2 - b^2 - c^2}{2bc} .
		\end{equation}
		Here, the integral exists as a Lebesgue integral when $\nu \in \intervaloo{1/2}{\infty}$, but as an improper integral when $\nu \in \intervaloc{-1/2}{1/2}$.
	\end{theorem}
In fact, Proposition \ref{prop:explicit formula for kernel} immediately follows from Theorem \ref{thm:explicit formula Legendre ver} by rewriting the Legendre functions as the Gegenbauer functions by Proposition \ref{prop:Gegenbauer and Legendre} and substituting $a = \rho$ and $b = c = s$.

	Next, we show that 
	\begin{equation} 
		\int_{r \in \intervaloo{0}{\infty}}f(r) \gmns(r) r^{2\nu+1} \, dr 
		= \int_{\rho \in \intervaloo{0}{\infty}} \Hankelnu{\nu} f(\rho) \Hankelnu{\nu}\gmns(\rho) \rho^{2\nu+1} \, d\rho 
	\end{equation}
	holds using Proposition \ref{prop:Parseval 3}.
	In order to apply Proposition \ref{prop:Parseval 3}, 
	it suffices to show that
	\begin{itemize}
		\item $\gmns \in \Lpnunu{1}{\nu+1/2}{\nu-3/2}$, 
		\item $r^{\nu+1/2} \gmns(r) \to 0$ as $r \to \infty$, 
		\item $\Hankelnu{\nu} \gmns \in \Lpnunu{1}{\nu+1/2}{2\nu+1}$.
	\end{itemize}
	The first and second properties are immediate from \eqref{eq:BesselJ upper bound}, and
	the third one is a consequence of Propositions \ref{prop:behavior of CD} and \ref{prop:explicit formula for kernel}.
	In fact, combining these propositions, we get
	\begin{align}
		\abs{ \Hankelnu{\nu} \gmns(\rho)  } 
		&\lesssim_{\mu, \nu, s} 
		\rho^{-1}
		\shortintertext{for $\rho \in (0, s)$, }
		\abs{ \Hankelnu{\nu} \gmns(\rho)  } 
		&\lesssim_{\mu, \nu, s} 
		\begin{dcases}
			(2s - \rho)^{(2\nu - 1)/2} , & \nu \in \intervaloo{-1/2}{1/2} \text{ or } \mu \in \N , \\
			\log{\mleft(\frac{2s}{2s - \rho}\mright)} , & \nu = 1/2 \text{ and } \mu \not \in \N , \\
			1  , & \nu \in \intervaloo{1/2}{\infty} \text{ and } \mu \not \in \N 
		\end{dcases}
		\shortintertext{for $\rho \in (s, 2s)$, }
		\abs{ \Hankelnu{\nu} \gmns(\rho)  } 
		&\lesssim_{\mu, \nu, s} 
		\begin{dcases}
			0 , & \mu \in \N , \\
			(\rho - 2s)^{(2\nu - 1)/2} , & \nu \in \intervaloo{-1/2}{1/2} \text{ and } \mu \not \in \N , \\
			\log{\mleft(\frac{2s}{\rho - 2s}\mright)} , & \nu = 1/2 \text{ and } \mu \not \in \N , \\
			1  , & \nu \in \intervaloo{1/2}{\infty} \text{ and } \mu \not \in \N 
		\end{dcases}
		\shortintertext{for $\rho \in (2s, 3s)$, and}
		\abs{ \Hankelnu{\nu} \gmns(\rho)  } 
		&\lesssim_{\mu, \nu, s} 
		\begin{cases}
			0 , & \mu \in \N , \\
			\rho^{ -2(\mu+\nu+1)} , &  \mu \not \in \N 
		\end{cases}
	\end{align}
	for $\rho \in (3s, \infty)$, 
	where $A \lesssim_{\mu, \nu, s} B$ means that there exists a constant $C_{\mu, \nu, s} > 0$ depending only on $\mu, \nu, s$ such that $A \leq C_{\mu, \nu, s}B$.
	As a consequence, we have $\Hankelnu{\nu} \gmns \in \Lpnunu{1}{\nu+1/2}{2\nu+1}$ whenever $\mu \in \intervalco{0}{\infty}$ and $\nu \in \intervaloo{-1/2}{\infty}$.
	\begin{remark}
		Our proof is substantially different from that of \citet{BSS2015} for the case $(\mu, \nu) = (k, d/2 - 1)$.
		Roughly speaking, their proof is as follows.
		Let $\varphi \colon \R \to \C$ be an even Schwartz function, and let $Y_{k, d} \colon \S^{d-1} \to \C$ be a spherical harmonic polynomial of degree $k$ in $d$ variables such that $\norm{Y_{k, d}}_{L^2(\S^{d-1})} = 1$.
		We consider the left-hand side of the smoothing estimate \eqref{eq:smoothing Schrodinger} with the initial data $u_0 \in L^2(\R^d)$ given by
		\begin{equation}
			\widehat{u_0}(r \theta) = r^{-(d-1)/2} \varphi(r) Y_{k, d}(\theta) , \quad (r, \theta) \in \intervaloo{0}{\infty} \times \S^{d-1}.  
		\end{equation}
		In the proof of Theorem \ref{thm:Walther}, \citet{Wal2002} showed that 
		\begin{equation}
			\int_{(x, t) \in \R^d \times \R}  w(\abs{x}) \abs{\psi(\abs{\nabla}) e^{it\Laplacian} u_0(x)}^2 \, dx \, dt
			= \int_0^\infty r^{-1} \psi(r)^2 \Tnu{k + d/2 - 1} w(r) \abs{\varphi(r)}^2 \, dr 
		\end{equation}
		holds using some identities involving the Fourier transform and the Bessel functions.
		On the other hand, \citet{BSS2015} calculated the same integral using a different method based on the so-called Funk--Hecke theorem, and showed that
		\begin{equation}
			\int_{(x, t) \in \R^d \times \R}  w(\abs{x}) \abs{\psi(\abs{\nabla}) e^{it\Laplacian} u_0(x)}^2 \, dx \, dt
			= \int_0^\infty r^{-1} \psi(r)^2 \Umunu{k}{d/2 - 1} \Hankelnu{d/2 - 1} w(r) \abs{\varphi(r)}^2 \, dr
		\end{equation}
		holds. Comparing these, we conclude that
		\begin{equation}
			\Tnu{k + d/2 - 1} w = \Umunu{k}{d/2 - 1} \Hankelnu{d/2 - 1} w 
		\end{equation}
		holds.
	\end{remark}
	\subsection{\texorpdfstring{Proofs of the inequalities \ref{item:inequality 2}, \ref{item:inequality 3}, \ref{item:inequality 4}, \ref{item:inequality 5}}{Proofs of the inequalities (II), (III), (IV), (V)}} \label{subsection:proofs of the inequalities}
	In order to prove \ref{item:inequality 2}, \ref{item:inequality 3}, \ref{item:inequality 4}, \ref{item:inequality 5},
	we first show the corresponding inequalities for the kernel $\Kernelmunu{\mu}{\nu}$. 
	For simplicity, we assume 
	\begin{equation}
		\Kernelmunu{\mu}{\nu}(r, s) = 0
	\end{equation}
	when $r = 2s$ in this section.
	\begin{proposition} \label{prop:ineq for kernel}
		The following hold:
		\begin{enumerate}[label=\textup{(\roman*)}, start=2]
			\item \label{item:inequality 2 kernel}
			Let $\nu \in \intervalco{1/2}{\infty}$. Then 
			\begin{equation}
				\intervaloo{0}{\infty} \ni s \longmapsto \Kernelmunu{0}{\nu}(r, s)
			\end{equation}
			is non-decreasing for each fixed $r \in \intervaloo{0}{\infty}$.
			Moreover, 
			\begin{equation}
				\set{ r \in \intervaloo{0}{\infty} }{ \Kernelmunu{0}{\nu}(r, s_2) - \Kernelmunu{0}{\nu}(r, s_1) > 0 } 
			\end{equation}
			has a positive Lebesgue measure for every $s_1, s_2 \in \intervaloo{0}{\infty}$ such that $s_1 < s_2$.
			\item \label{item:inequality 3 kernel}
			Let $\nu \in \intervaloo{-1/2}{\infty}$. Then 
			\begin{equation}
				\intervaloo{0}{\infty} \ni s \longmapsto \Kernelmunu{0}{\nu}(r, s) + \Kernelmunu{1}{\nu}(r, s)
			\end{equation}
			is non-decreasing for each fixed $r \in \intervaloo{0}{\infty}$.
			Moreover, 
			\begin{equation}
				\set{ r \in \intervaloo{0}{\infty} }{ \Kernelmunu{0}{\nu}(r, s_2) + \Kernelmunu{1}{\nu}(r, s_2) - ( \Kernelmunu{0}{\nu}(r, s_1) + \Kernelmunu{1}{\nu}(r, s_1) ) > 0 } 
			\end{equation}
			has a positive Lebesgue measure for every $s_1, s_2 \in \intervaloo{0}{\infty}$ such that $s_1 < s_2$.
			\item \label{item:inequality 4 kernel}
			Let $\mu \in \N_{\geq 1}$ and $\nu \in \intervalco{0}{\infty}$. 
			Then we have
			\begin{equation}
				\Kernelmunu{\mu}{\nu}(r, s) \leq \Kernelmunu{0}{\nu}(r, s) 
			\end{equation}
			for every $r, s \in \intervaloo{0}{\infty}$. Moreover, 
			\begin{equation}
				\set{ r \in \intervaloo{0}{\infty} }{ \Kernelmunu{0}{\nu}(r, s) - \Kernelmunu{\mu}{\nu}(r, s) > 0 } 
			\end{equation}
			has a positive Lebesgue measure for every $s \in \intervaloo{0}{\infty}$.
			\item \label{item:inequality 5 kernel}
			Let $\mu \in \intervaloc{0}{1}$ and $\nu \in \intervaloo{-1/2}{\infty}$ be such that $\mu + 2\nu \geq 0$. 
			Then we have
			\begin{equation}
				\Kernelmunu{\mu}{\nu}(r, s) \leq \Kernelmunu{0}{\nu}(r, s) 
			\end{equation}
			for every $r, s \in \intervaloo{0}{\infty}$. Moreover, 
			\begin{equation}
				\set{ r \in \intervaloo{0}{\infty} }{ \Kernelmunu{0}{\nu}(r, s) - \Kernelmunu{\mu}{\nu}(r, s) > 0 } 
			\end{equation}
			has a positive Lebesgue measure for every $s \in \intervaloo{0}{\infty}$.
		\end{enumerate}
	\end{proposition}
	\begin{proof}[Proof of Proposition \ref{prop:ineq for kernel}]
		Recall that we have
		\begin{equation}
			\Cmunu{0}{\nu}(x) \equiv 1, \quad \Cmunu{1}{\nu}(x) \equiv x , 
		\end{equation}
		so that
		\begin{align}
			\Kernelmunu{0}{\nu}(r, s) &= \dfrac{\pi^{1/2} r^{2\nu}}{2^{\nu} \Gamma(\nu+1/2) }
			\times
			\begin{cases}
				 (1 - r^2/(4s^2))^{(2\nu - 1)/2}  , & 0 < r < 2s , \\
				0 , & r \geq 2s > 0 ,
			\end{cases}  \\
			\Kernelmunu{0}{\nu}(r, s) + \Kernelmunu{1}{\nu}(r, s) &= \dfrac{2 \pi^{1/2} r^{2\nu}}{2^{\nu} \Gamma(\nu+1/2) }
			\times
			\begin{cases}
			 (1 - r^2/(4s^2))^{(2\nu + 1)/2}  , & 0 < r < 2s , \\
				0 , & r \geq 2s > 0 .
			\end{cases}
		\end{align}
		The monotonicities in \ref{item:inequality 2 kernel}, \ref{item:inequality 3 kernel} are immediate from these.
		It is also easy to see that
		\begin{gather}
			\set{ r \in \intervaloo{0}{\infty} }{ \Kernelmunu{0}{\nu}(r, s_2) - \Kernelmunu{0}{\nu}(r, s_1) > 0 } 
			= \begin{cases}
				\intervalco{2s_1}{2s_2} , & \nu = 1/2 , \\
				\intervaloo{0}{2s_2} , & \nu \in \intervaloo{1/2}{\infty} ,
			\end{cases} 
			\label{eq:identity 3 kernel set} \noeqref{eq:identity 3 kernel set}\\
			\set{ r \in \intervaloo{0}{\infty} }{ \Kernelmunu{0}{\nu}(r, s_2) + \Kernelmunu{1}{\nu}(r, s_2) - ( \Kernelmunu{0}{\nu}(r, s_1) + \Kernelmunu{1}{\nu}(r, s_1) ) > 0 } 
			= \intervaloo{0}{2s_2} 
			\label{eq:identity 4 kernel set} \noeqref{eq:identity 4 kernel set}
		\end{gather}
		hold for every $s_1, s_2 \in \intervaloo{0}{\infty}$ such that $s_1 < s_2$.
		
		Next, we prove \ref{item:inequality 4 kernel}. When $\mu = k \in \N_{\geq 1}$, we have
		\begin{equation} 
			\Kernelmunu{k}{\nu}(r, s) = \dfrac{\pi^{1/2} r^{2\nu}}{2^{\nu} \Gamma(\nu+1/2) }
			\times
			\begin{cases}
				(1 - r^2/(4s^2))^{(2\nu - 1)/2} \Cmunu{k}{\nu}(1 - r^2/(2s^2)) , & 0 < r < 2s , \\
				0 , & r \geq 2s > 0 .
			\end{cases}  
		\end{equation}
		Therefore, the desired inequality follows from \eqref{eq:Cmunu upper bound integer}.
		We also have 
		\begin{equation} 	\label{eq:identity 2a kernel set} \noeqref{eq:identity 2a kernel set}
			\set{ r \in \intervaloo{0}{\infty} }{ \Kernelmunu{0}{\nu}(r, s) - \Kernelmunu{k}{\nu}(r, s) > 0 } 
			= \intervaloo{0}{2s} \setminus \Nullset_{k, s} ,
		\end{equation}
		where
		\begin{equation}
			\Nullset_{k, s} \coloneqq \set{ r \in \intervaloo{0}{2s} }{ \Cmunu{k}{\nu}(1 - r^2/(2s^2)) = 1 } ,
		\end{equation}
		which contains at most $k$ points.
		
		Finally, we prove \ref{item:inequality 5 kernel}.
		Observe that the desired inequality and
		\begin{equation} \label{eq:identity 2b kernel set} \noeqref{eq:identity 2b kernel set}
			\set{ r \in \intervaloo{0}{\infty} }{ \Kernelmunu{0}{\nu}(r, s) - \Kernelmunu{\mu}{\nu}(r, s) > 0 } 
			= \begin{cases}
				\intervaloo{2s}{\infty}, & \mu + 2\nu = 0 , \\
				\intervaloo{0}{\infty} \setminus \{2s\}, & \mu + 2\nu > 0 , \, \mu \ne 1 , \\
				\intervaloo{0}{2s} , & \mu = 1 ,
			\end{cases}
		\end{equation}
	immediately follow from \eqref{eq:Cmunu upper bound} and \eqref{eq:Dmunu positive}, which completes the proof.
	\end{proof}	
	Now \ref{item:inequality 2}, \ref{item:inequality 3}, \ref{item:inequality 4}, \ref{item:inequality 5} are easily derived from \ref{item:inequality 2 kernel}, \ref{item:inequality 3 kernel}, \ref{item:inequality 4 kernel}, \ref{item:inequality 5 kernel}, respectively. 
	We shall only prove \ref{item:inequality 2}, since the proofs of the others are similar.
	\begin{proof}[Proof of \ref{item:inequality 2}]
		Let $\nu \in \intervalco{1/2}{\infty}$, and let
		$f \in \Lpnunu{1}{2\nu+1}{0}$ be such that $\Hankelnu{\nu} f$ is non-negative, that is, 
		\begin{equation}
			\Hankelnu{\nu} f_\varepsilon(r) \geq 0
		\end{equation}
		holds for every $r \in \intervaloo{0}{\infty}$, where $f_\varepsilon$ is given by
		\begin{equation}
			f_\varepsilon(r) \coloneqq f(r) \exp(- \varepsilon r^2 ) 
		\end{equation}
		for each $\varepsilon \in \intervaloo{0}{\infty}$. 
		Now fix $s_1, s_2 \in \intervaloo{0}{\infty}$ satisfying $s_1 < s_2$ arbitrarily.
		Then
		we have
		\begin{equation}
			\Tnu{\nu} f_\varepsilon(s_2) - \Tnu{\nu} f_\varepsilon(s_1) = \Umunu{0}{\nu} \Hankelnu{\nu} f_\varepsilon(s_2) - \Umunu{0}{\nu} \Hankelnu{\nu} f_\varepsilon(s_1)
		\end{equation}
		by \ref{item:identity 1}, and 
		\begin{equation}
			\lim_{\varepsilon \downarrow 0} {(\Tnu{\nu} f_\varepsilon(s_2) - \Tnu{\nu} f_\varepsilon(s_1))} = \Tnu{\nu}f(s_2) - \Tnu{\nu}f(s_1)
		\end{equation}
		by the dominated convergence theorem.
		On the other hand, combining \ref{item:inequality 2 kernel} and the non-negativity of $\Hankelnu{\nu} f$, it is clear that
		\begin{equation}
			\Umunu{0}{\nu} \Hankelnu{\nu} f_\varepsilon(s_2) - \Umunu{0}{\nu} \Hankelnu{\nu} f_\varepsilon(s_1) \geq 0
		\end{equation} 
		holds for every $\varepsilon \in \intervaloo{0}{\infty}$.
		Thus, we conclude that 
		\begin{equation}
			\Tnu{\nu}f(s_2) - \Tnu{\nu}f(s_1) \geq 0
		\end{equation}
		holds. Furthermore, when $\Hankelnu{\nu} f$ is strictly positive, then we also have
		\begin{equation}
			F(r) \coloneqq \liminf_{\varepsilon \downarrow 0} \Hankelnu{\nu} f_\varepsilon(r) > 0
		\end{equation}
		for almost every $r > 0$. Hence, using Fatou's lemma, we get
		\begin{equation}
			\liminf_{\varepsilon \downarrow 0} {(\Umunu{0}{\nu} \Hankelnu{\nu} f_\varepsilon(s_2) - \Umunu{0}{\nu} \Hankelnu{\nu} f_\varepsilon(s_1) )} \geq \Umunu{0}{\nu} F(s_2) - \Umunu{0}{\nu} F(s_1) > 0 , 
		\end{equation}
		so that 
		\begin{equation}
			\Tnu{\nu}f(s_2) - \Tnu{\nu}f(s_1) \geq \Umunu{0}{\nu} F(s_2) - \Umunu{0}{\nu} F(s_1) > 0
		\end{equation}
		holds. This completes the proof. 
	\end{proof}	
	
	\section{\texorpdfstring{Proof of Theorem \ref{main thm:Tnuf for completely monotone f}}{Proof of Theorem 1.6}} \label{section:completely monotone}
	In this section, we prove Theorem \ref{main thm:Tnuf for completely monotone f}. 
	Throughout this section, we write
	\begin{equation}
		\varphi(r) \coloneqq \exp{(- r^2 / 2)} 
	\end{equation}
	for simplicity.
	First, we consider the case $f = \varphi$.
	In this case, it is well known that
	\begin{equation}
		\Hankelnu{\nu} \varphi(r) = \varphi(r) > 0 
		\label{eq:Hankelnu Gaussian}
	\end{equation}
	holds (see \cite[\href{https://dlmf.nist.gov/10.22.E51}{10.22.51}]{DLMF}, \cite[6.631.4]{GR2014}). 
	Therefore, using Theorem \ref{main thm:identity and inequalities}, we obtain the identity
	\begin{equation}
		\Tnu{\nu} \varphi(s) = \Umunu{0}{\nu} \varphi(s) 
	\end{equation}
	and the strict inequalities corresponding to \ref{item:inequality 2}, \ref{item:inequality 3}, \ref{item:inequality 4}, \ref{item:inequality 5}.
	Now notice that the kernel $\Kernelmunu{0}{\nu}$ satisfies
	\begin{equation} \label{eq:kernel homogeneous}
		\Kernelmunu{0}{\nu}(ar, as) = a^{2\nu} \Kernelmunu{0}{\nu}(r, s)
	\end{equation}
	for every $a, r, s \in \intervaloo{0}{\infty}$.
	Therefore, by changing variable of integration, we get
	\begin{align}
		\Umunu{0}{\nu} \varphi(s)
		&= \int_{r \in \intervaloo{0}{\infty}} \Kernelmunu{0}{\nu}(r, s) \varphi(r) \, dr \\
		\underset{\eqref{eq:kernel homogeneous}}&= s^{2\nu+1} \int_{r \in \intervaloo{0}{\infty}} \Kernelmunu{0}{\nu}(r, 1) \varphi(rs) \, dr ,
	\end{align}
	so that 
	\begin{equation}
		s^{-\nu-1/2} \Tnu{\nu} \varphi(s^{1/2}) = s^{-\nu-1/2} \Umunu{0}{\nu} \varphi(s^{1/2}) = \int_{r \in \intervaloo{0}{\infty}} \Kernelmunu{0}{\nu}(r, 1) \exp{(-r^2 s / 2)} \, dr
	\end{equation}
	holds. Since $\Kernelmunu{0}{\nu}$ is non-negative, we conclude that 
	\begin{equation}
		s \longmapsto s^{-\nu-1/2} \Tnu{\nu} \varphi(s^{1/2})
	\end{equation}
	is completely monotone by Theorem \ref{thm:completely monotone functions}.
	Moreover, actually it is known that
	\begin{equation}
		\Tnu{\nu} \varphi(s) = \pi s \exp{(-s^2)} \BesselI{\nu}(s^2) 
		\label{eq:Tnu Gaussian}
	\end{equation}
	holds for every $\nu \in \intervalco{-1/2}{\infty}$ (see \cite[\href{https://dlmf.nist.gov/10.22.E67}{10.22.67}]{DLMF}, \cite[6.633.2]{GR2014}).
	In summary, we have:
	\begin{example} \label{example:Tnu Gaussian}
		The following hold.
		\begin{enumerate}[label=\textup{(\arabic*a)}]
			\item \label{item:Tnu Gaussian 1}
			For each fixed $\nu \in \intervaloo{-1/2}{\infty}$,
			\begin{equation}
				s \longmapsto s^{-\nu} \exp{(-s)} \BesselI{\nu}(s)
			\end{equation}
			is completely monotone on $\intervaloo{0}{\infty}$.
			\item \label{item:Tnu Gaussian 2}
			For each fixed $\nu \in \intervalco{1/2}{\infty}$, 
			\begin{equation}
				s \longmapsto s^{1/2} \exp{(-s)} \BesselI{\nu}(s)
			\end{equation}
			is strictly increasing on $\intervaloo{0}{\infty}$.
			\item \label{item:Tnu Gaussian 3}
			For each fixed $\nu \in \intervaloo{-1/2}{\infty}$, 
			\begin{equation}
				s \longmapsto s^{1/2} \exp{(-s)} ( \BesselI{\nu}(s) + \BesselI{\nu+1}(s) )
			\end{equation}
			is strictly increasing on $\intervaloo{0}{\infty}$.
			\item \label{item:Tnu Gaussian 4}
			For each fixed $s \in \intervaloo{0}{\infty}$, 
			\begin{equation}
				\nu \longmapsto \BesselI{\nu}(s)
			\end{equation}
			is strictly decreasing on $\intervalco{0}{\infty}$.
			\item \label{item:Tnu Gaussian 5}
			For each fixed $\nu \in \intervaloo{-1/2}{0}$ and $s \in \intervaloo{0}{\infty}$, 
			\begin{equation}
				\BesselI{-\nu}(s) < \BesselI{\nu}(s)
			\end{equation}
			holds.
		\end{enumerate}
	\end{example}
	Now we prove Theorem \ref{main thm:Tnuf for completely monotone f}.
	\begin{proof}[Proof of Theorem \ref{main thm:Tnuf for completely monotone f}]
		For simplicity, we write
		\begin{equation}
			\varphi_a(r) \coloneqq \varphi(ar) = \exp{(-a^2 r^2 / 2)} 
		\end{equation}
		for each $a \in \intervalco{0}{\infty}$.
		Let $f \in \Lpnunu{1}{2\nu_0+1}{0}$ be such that $r \longmapsto f(r^{1/2})$ is completely monotone and not identically zero. 
		By Theorem \ref{thm:completely monotone functions}, there exists a Borel measure $\lambda$ on $\intervalco{0}{\infty}$ such that
		\begin{equation}
			f(r) = \int_{a \in \intervalco{0}{\infty}} \varphi_a(r) \, d\lambda(a) 
		\end{equation}
		holds for every $r \in \intervaloo{0}{\infty}$.
		Since $f$ is not identically zero, this measure satisfies $\lambda(\intervalco{0}{\infty}) > 0$. 
		Furthermore, the assumption $f \in \Lpnunu{1}{2\nu_0+1}{0}$ implies $\lambda(\{0\}) = 0$, so that we have
		\begin{equation} \label{eq:f by phi}
			f(r) = \int_{a \in \intervaloo{0}{\infty}} \varphi_a(r) \, d\lambda(a) .
		\end{equation}
		Hence, changing order of integration by Tonelli's theorem, we get
		\begin{equation} \label{eq:Tnuf by Tnuphi}
			\Tnu{\nu}f(s) 
			= \int_{a \in \intervaloo{0}{\infty}} \Tnu{\nu} \varphi_a(s) \, d \lambda(a) 
		\end{equation}
		for each $\nu \in \intervalco{\nu_0}{\infty}$.
		Now note that
		\begin{equation}
			\Hankelnu{\nu} \varphi_a = a^{-2(\nu+1)} \varphi_{1/a} 
		\end{equation}
		holds for each $a \in \intervaloo{0}{\infty}$, which implies
		\begin{align}
			\Tnu{\nu} \varphi_a(s) 
			\underset{\ref{item:identity 1}}&= a^{-2(\nu+1)} \Umunu{0}{\nu} \varphi_{1/a}(s) \\
			&= a^{-2(\nu+1)} \int_{r \in \intervaloo{0}{\infty}} \Kernelmunu{0}{\nu}(r, s) \varphi_{1/a}(r) \, dr \\
			\underset{\eqref{eq:kernel homogeneous}}&= a^{-1} s^{2\nu+1} \int_{r \in \intervaloo{0}{\infty}} \Kernelmunu{0}{\nu}(r, 1/a) \varphi(rs) \, dr .
			\label{eq:Tnuphi}
		\end{align}
		Therefore, substituting \eqref{eq:Tnuphi} into \eqref{eq:Tnuf by Tnuphi} and changing order of integration again, we obtain
		\begin{equation}
			\Tnu{\nu}f(s) 
			= s^{2\nu+1} \int_{r \in \intervaloo{0}{\infty}} \mleft(\int_{a \in \intervaloo{0}{\infty}} a^{-1} \Kernelmunu{0}{\nu}(r, 1/a) \, d \lambda(a) \mright) \varphi(rs)  \, dr , 
		\end{equation}
		or equivalently
		\begin{equation}
			s^{-(\nu+1/2)} \Tnu{\nu}f(s^{1/2}) 
			= \int_{r \in \intervaloo{0}{\infty}} \mleft(\int_{a \in \intervaloo{0}{\infty}} a^{-1} \Kernelmunu{0}{\nu}(r, 1/a) \, d \lambda(a) \mright) \exp{(-r^2 s / 2)}  \, dr .
		\end{equation}
		Thus, by Theorem \ref{thm:completely monotone functions}, we conclude that
		\begin{equation}
			s \longmapsto s^{-(\nu+1/2)} \Tnu{\nu}f(s^{1/2})
		\end{equation}
		is completely monotone. 
		This shows \ref{item:Tnuf for completely monotone 1}. 
		
		In order to prove the inequalities \ref{item:Tnuf for completely monotone 2}, \ref{item:Tnuf for completely monotone 3}, \ref{item:Tnuf for completely monotone 4}, \ref{item:Tnuf for completely monotone 5}, one can use Theorem \ref{main thm:identity and inequalities} by showing that $\Hankelnu{\nu} f$ is strictly positive.
		Alternatively, they follow from 
		\begin{equation}
			\Tnu{\nu}f(s) 
			\underset{\eqref{eq:Tnuf by Tnuphi}}= \int_{a \in \intervaloo{0}{\infty}} \Tnu{\nu} \varphi_a(s) \, d \lambda(a)
			\underset{\eqref{eq:Tnu Gaussian}}= \pi s \int_{a \in \intervaloo{0}{\infty}} a^{-2} \exp{(-s^2/a^2)} \BesselI{\nu}(s^2/a^2) \, d \lambda(a) 
		\end{equation}
		and \ref{item:Tnu Gaussian 2}, \ref{item:Tnu Gaussian 3}, \ref{item:Tnu Gaussian 4}, \ref{item:Tnu Gaussian 5} of Example \ref{example:Tnu Gaussian}.
		This completes the proof of Theorem \ref{main thm:Tnuf for completely monotone f}.
	\end{proof}
	As an application of Theorem \ref{main thm:Tnuf for completely monotone f}, for example, we can show the following: 
	\begin{example} \label{example:Tnu inhomogeneous power}
		The following hold.
		\begin{enumerate}[label=\textup{(\arabic*b)}]
			\item \label{item:Tnu inhomogeneous power 1}
			For each fixed $\nu \in \intervaloo{-1/2}{\infty}$,
			\begin{equation}
				s \longmapsto s^{-\nu} \BesselI{\nu}(s^{1/2}) \BesselK{\nu}(s^{1/2})
			\end{equation}
			is completely monotone on $\intervaloo{0}{\infty}$.
			\item \label{item:Tnu inhomogeneous power 2}
			For each fixed $\nu \in \intervalco{1/2}{\infty}$, 
			\begin{equation}
				s \longmapsto s \BesselI{\nu}(s) \BesselK{\nu}(s)
			\end{equation}
			is strictly increasing on $\intervaloo{0}{\infty}$.
			\item \label{item:Tnu inhomogeneous power 3}
			For each fixed $\nu \in \intervaloo{-1/2}{\infty}$, 
			\begin{equation}
				s \longmapsto s ( \BesselI{\nu}(s) \BesselK{\nu}(s) + \BesselI{\nu+1}(s) \BesselK{\nu+1}(s) )
			\end{equation}
			is strictly increasing on $\intervaloo{0}{\infty}$.
			\item \label{item:Tnu inhomogeneous power 4}
			For each fixed $s \in \intervaloo{0}{\infty}$, 
			\begin{equation}
				\nu \longmapsto \BesselI{\nu}(s) \BesselK{\nu}(s)
			\end{equation}
			is strictly decreasing on $\intervalco{0}{\infty}$.
			\item \label{item:Tnu inhomogeneous power 5}
			For each fixed $\nu \in \intervaloo{-1/2}{0}$ and $s \in \intervaloo{0}{\infty}$, 
			\begin{equation}
				\BesselI{-\nu}(s) \BesselK{-\nu}(s) < \BesselI{\nu}(s) \BesselK{\nu}(s)
			\end{equation}
			holds.
		\end{enumerate}
	\end{example}
	Example \ref{example:Tnu inhomogeneous power} follows from Theorem \ref{main thm:Tnuf for completely monotone f} by letting 
	\begin{equation}
		f(r) = (1+r^2)^{-1} , 
	\end{equation}
	since 
	\begin{equation}
		f(r) = (1+r^2)^{-1} \longmapsto	\Tnu{\nu}f(s) = \pi s \BesselI{\nu}(s) \BesselK{\nu}(s) 
	\end{equation}
	holds for every $\nu \in \intervalco{-1/2}{\infty}$ (see \cite[6.535]{GR2014}).
	In fact, it is clear that $f \in \Lpnunu{1}{0}{0} = L^1$ and that $r \mapsto f(r^{1/2}) = (1+r)^{-1}$ is completely monotone.
	\begin{remark}
		Examples \ref{example:Tnu Gaussian} and \ref{example:Tnu inhomogeneous power} themselves are already known, even though our proofs are new. Here we give references.
		\begin{itemize}
			\item \ref{item:Tnu Gaussian 1} was observed by \citet[\textsection 2]{Nas1978}. 
			\item \citet[p.\ 587]{Bar2010-2} observed an analogue of \ref{item:Tnu Gaussian 1} for the modified Bessel function of the second kind $\BesselK{\nu}$, that is, 
			\begin{equation}
				s \longmapsto s^{-\nu} \exp{(s)} \BesselK{\nu}(s)
			\end{equation}
			is completely monotone for each $\nu \in \intervaloo{-1/2}{\infty}$.
			Since a product of two completely monotone functions is also completely monotone, we see that
			\begin{equation}
				s \longmapsto s^{-\nu} \exp{(-s)} \BesselI{\nu}(s) \times s^{-\nu} \exp{(s)} \BesselK{\nu}(s) = s^{-2 \nu} \BesselI{\nu}(s) \BesselK{\nu}(s)
			\end{equation}
			is completely monotone.
			Note that this implies \ref{item:Tnu inhomogeneous power 1}, since $s \longmapsto f(s^{1/2})$ is also completely monotone when $f$ is (see \cite[Theorem E']{LN1983}, \cite[Theorem 2]{MS2001}, \cite[Theorem 2.5, Corollary 2.6]{Sza2026}, for example).				
			\item
			\citet[(2.6)]{Bar2010-2} pointed out that \ref{item:Tnu Gaussian 2} can be derived from \citeauthor{Gro1932}'s result \cite[(5)]{Gro1932}
			\begin{equation}
				\frac{d}{ds} \BesselI{\nu}(s) > (1 - 1/(2s)) \BesselI{\nu}(s) ,
				\quad s \in \intervaloo{0}{\infty}, \nu \in \intervalco{1/2}{\infty} , 
			\end{equation}
			since
			\begin{equation}
				\frac{d}{ds} (s^{1/2} \exp{(-s)} \BesselI{\nu}(s)) = s^{1/2} \exp{(-s)} \mleft( \frac{d}{ds} \BesselI{\nu}(s) - (1 - 1/(2s)) \BesselI{\nu}(s) \mright) .
			\end{equation}
			\item 
			\ref{item:Tnu inhomogeneous power 2} was observed by \citet[Theorem 4.1]{Har1977}. 
			\citet[Theorem 4.2]{Har1977} also showed that it is not monotone when $\nu \in \intervalco{0}{1/2}$.
			More precisely, for each $\nu \in \intervalco{0}{1/2}$, there exists $s_\nu \in \intervaloo{0}{\infty}$ such that 
			\begin{equation}
				s \longmapsto s \BesselI{\nu}(s) \BesselK{\nu}(s)
			\end{equation}
			is strictly increasing on $\intervaloc{0}{s_\nu}$ and strictly decreasing on $\intervalco{s_\nu}{\infty}$.
			\item For \ref{item:Tnu Gaussian 3} and \ref{item:Tnu inhomogeneous power 3}, see Proposition \ref{prop:square sum of BesselJ}.
			\item \ref{item:Tnu Gaussian 4} was observed by \citet{Coc1967, Reu1968, Jon1968} independently around the same time.
			Their works are inspired by \citet{Son1965}, which shows 
			\begin{equation}
				\BesselI{\nu + 1}(s) < \BesselI{\nu}(s)
			\end{equation}
			holds for every $\nu \in \intervaloo{-1/2}{\infty}$ and $s \in \intervaloo{0}{\infty}$.
			\item \ref{item:Tnu inhomogeneous power 4} was observed by \citet[Theorem 2]{BP2013}. 
			\item We can show \ref{item:Tnu Gaussian 5} and \ref{item:Tnu inhomogeneous power 5} easily and directly using \cite[\href{https://dlmf.nist.gov/10.27.E2}{10.27.2}, \href{https://dlmf.nist.gov/10.27.E3}{10.27.3}]{DLMF}
			\begin{align}
				\BesselI{\nu}(s) - \BesselI{-\nu}(s) &= - \frac{2}{\pi} \sin{(\nu \pi)} \BesselK{\nu}(s) , \\
				\BesselK{\nu}(s) &= \BesselK{-\nu}(s) .
			\end{align}
		\end{itemize}
	\end{remark}
	\section{Remarks} \label{section:remarks}
	\subsection{On the case \texorpdfstring{$\nu = -1/2$}{ν=-1/2}}\label{subsection:v=-1/2}
	For each $\mu \in \{ 0, 1 \}$, we define a linear operator $\Umunu{\mu}{-1/2}$ by
	\begin{equation}
		\Umunu{\mu}{-1/2} f(s) 
		\coloneqq (\pi / 2)^{1/2} \times \begin{cases}
			f(0) + f(2s) , & \mu = 0 , \\
			f(0) - f(2s) , & \mu = 1 .
		\end{cases}
	\end{equation}
	Then we have the following analogues of Theorems \ref{main thm:identity and inequalities} and \ref{main thm:Tnuf for completely monotone f}.
	\begin{proposition} \label{prop:identity and inequalities v=-1/2}
		Let $f \in \Lpnunu{1}{0}{0} = L^1$. Then the following hold.
		\begin{enumerate}[label=\textup{(\Roman**)}]
			\item \label{item:identity 1 v=-1/2}
			For each $\mu \in \{ 0, 1 \}$ and every $s \in \intervaloo{0}{\infty}$, 
			\begin{equation}
				\Tnu{\mu - 1/2} f(s) = \Umunu{\mu}{-1/2} \Hankelnu{-1/2} f(s)
			\end{equation}
			holds.
			\setcounter{enumi}{2}
			\item \label{item:inequality 3 v=-1/2}
			For every $s \in \intervaloo{0}{\infty}$, 
			\begin{equation}
				\Tnu{-1/2} f(s) + \Tnu{1/2} f(s) = 2 \int_0^\infty f(r) \, dr
			\end{equation}
			holds.
			\setcounter{enumi}{4}
			\item \label{item:inequality 5 v=-1/2}
			For each fixed $s \in \intervaloo{0}{\infty}$, 
			\begin{equation}
				\Tnu{1/2} f(s) \leq \Tnu{-1/2} f(s)
			\end{equation}
			holds if $\Hankelnu{-1/2} f(2s) \geq 0$, and the inequality is strict if $\Hankelnu{-1/2} f(2s) > 0$.
		\end{enumerate}
		\begin{enumerate}[label=\textup{(\arabic**)}]
			\item \label{item:Tnuf for completely monotone 1 v=-1/2}
			If $r \longmapsto f(r^{1/2})$ is completely monotone, then $s \longmapsto \Tnu{-1/2}f(s^{1/2})$ is also completely monotone.
		\end{enumerate}
	\end{proposition}
	\ref{item:identity 1 v=-1/2}, \ref{item:inequality 3 v=-1/2}, \ref{item:inequality 5 v=-1/2} are analogues of \ref{item:identity 1}, \ref{item:inequality 3}, \ref{item:inequality 5} of Theorem \ref{main thm:identity and inequalities}, and \ref{item:Tnuf for completely monotone 1 v=-1/2} is that of \ref{item:Tnuf for completely monotone 1} of Theorem \ref{main thm:Tnuf for completely monotone f}, respectively.
	\ref{item:identity 1 v=-1/2} follows from
	\begin{align}
		\pi r \BesselJ{-1/2}(r)^2 &= 2\cos^2{(r)} = 1 + \cos{(2r)} , \\
		\pi r \BesselJ{1/2}(r)^2 &= 2 \sin^2{(r)} = 1 - \cos{(2r)} , \\
		r^{1/2} \BesselJ{-1/2}(r) &= (2/\pi)^{1/2} \cos{(r)} ,
	\end{align}
	and \ref{item:inequality 3 v=-1/2}, \ref{item:inequality 5 v=-1/2}, \ref{item:Tnuf for completely monotone 1 v=-1/2} follow from \ref{item:identity 1 v=-1/2}. 
	\subsection{\texorpdfstring{An analogue of Theorem \ref{main thm:d+1 vs d} for the Dirac equations}{An analogue of Theorem 1.3 for the Dirac equations}} \label{subsection:Dirac}
	Using a similar argument, we can show an analogue of Theorem \ref{main thm:d+1 vs d} for the Dirac equations.
	Let $\Diracconstwpmd{w}{\psi}{m}{d}$ be the optimal constant of the inequality
	\begin{equation} 
		\int_{(x, t) \in \R^d \times \R} w(\abs{x}) \abs{ ( - \Laplacian  +  m^2 )^{-1/4} \psi(\abs{\nabla}) e^{- it \Diracopm{m} } u_0(x)}^2 \, dx \, dt \leq C \norm{u_0}^2_{L^2} ,
		\label{eq:smoothing Dirac}
	\end{equation}
	where $\Diracopm{m}$ denotes the Dirac operator with mass $m \geq 0$.
	\citet{Suz2025} showed the following.
	\begin{theorem}[{\cite[Theorem 1.6]{Suz2025}}] \label{thm:Walther Dirac}
		Let $d \geq 2$. Then we have
		\begin{equation}
			\Diracconstwpmd{w}{\psi}{m}{d} = \sup_{k \in \N} \sup_{s > 0} s^{-1} \psi(s)^2 \DiracTnum{k + d/2 - 1}{m}w(s) ,
		\end{equation}
		where $\DiracTnum{\nu}{m}$ is a sublinear operator defined by
		\begin{equation}
			\DiracTnum{\nu}{m}f(s) \coloneqq \Tnu{\nu} f(s) + \Tnu{\nu+1} f(s) + \dfrac{m}{ \sqrt{s^2 + m^2} } \abs{ \Tnu{\nu} f(s) - \Tnu{\nu+1} f(s) } .
		\end{equation}
	\end{theorem}
	Combining Theorems \ref{main thm:identity and inequalities} and \ref{thm:Walther Dirac}, we get the following.
	\begin{theorem} \label{thm:d+1 vs d Dirac}
		Let $d \geq 2$.
		Assume that all of the Fourier transforms of a spatial weight $w \in \Lpnunu{1}{d-1}{0}$ as radial functions on $\R^d$, $\R^{d+1}$, $\R^{d+2}$, $\R^{d+3}$ (i.e., $\Hankelnu{d/2-1} w$, $\Hankelnu{d/2-1/2} w$, $\Hankelnu{d/2} w$, $\Hankelnu{d/2+1/2} w$) are non-negative.
		Then we have
		\begin{equation}
			\Diracconstwpmd{w}{\psi}{m}{d+1} \leq \Diracconstwpmd{w}{\psi}{m}{d}
		\end{equation}
		for every smoothing function $\psi \colon \intervaloo{0}{\infty} \to \intervalco{0}{\infty}$.
	\end{theorem}
	\begin{proof}[Proof of Theorem \ref{thm:d+1 vs d Dirac}]
		Since $\Hankelnu{d/2 - 1} w$ and $\Hankelnu{d/2} w$ are non-negative by the assumption, 
		\ref{item:inequality 4} and \ref{item:inequality 5} of Theorem \ref{main thm:identity and inequalities} imply that
		\begin{align}
			\Tnu{\mu + d/2 - 1} w(s) &\leq \Tnu{d/2 - 1} w(s) , \\
			\Tnu{\mu + d/2} w(s) &\leq \Tnu{d/2} w(s) ,
		\end{align}
		and consequently
		\begin{equation} \label{eq:Tnuw ineq Dirac}
			\DiracTnum{\mu + d/2 - 1}{m} w(s) \leq \DiracTnum{d/2 - 1}{m} w(s)
		\end{equation}
		holds for every $\mu \in \intervaloo{0}{1} \cup \N$ and $s \in \intervaloo{0}{\infty}$.
		Therefore, Theorem \ref{thm:Walther Dirac} gives us
		\begin{equation}
			\Diracconstwpmd{w}{\psi}{m}{d}
			\underset{\text{Theorem \ref{thm:Walther Dirac}}}= \sup_{k \in \N} \sup_{s > 0} s^{-1} \psi(s)^2 \DiracTnum{k + d/2 - 1}{m}w(s) 
			\underset{\eqref{eq:Tnuw ineq Dirac}}= \sup_{s > 0} s^{-1} \psi(s)^2 \DiracTnum{d/2 - 1}{m}w(s) .
		\end{equation}
		By the same argument using the non-negativity of $\Hankelnu{d/2 - 1/2} w$ and $\Hankelnu{d/2 + 1/2} w$, we also get
		\begin{equation}
			\Diracconstwpmd{w}{\psi}{m}{d+1} = \sup_{s > 0} s^{-1} \psi(s)^2 \DiracTnum{d/2 - 1/2}{m}w(s) .
		\end{equation}
		Now we use \eqref{eq:Tnuw ineq Dirac} with $\mu = 1/2$ and conclude that  
		\begin{equation}
			\Diracconstwpmd{w}{\psi}{m}{d+1} \leq \Diracconstwpmd{w}{\psi}{m}{d} 
		\end{equation}
		holds, as desired.
	\end{proof}
	As in Theorem \ref{main thm:d+1 vs d}, the assumption of Theorem \ref{thm:d+1 vs d Dirac} is satisfied in the cases \eqref{eq:type A}, \eqref{eq:type B}, \eqref{eq:type C}.
	See \cite[Theorems 1.7, 1.8, 1.9, 1.10]{Suz2025} for the explicit values of $\Diracconstwpmd{w}{\psi}{m}{d}$ in these cases.
	
	\subsection{\texorpdfstring{\ref{item:inequality 5} implies \ref{item:inequality 3}}{(V) implies (III)}} \label{subsection:differentiation}
	The following Proposition \ref{prop:square sum of BesselJ} reveals that \ref{item:inequality 3} can also be derived from \ref{item:inequality 5} rather than \ref{item:identity 1}.
	\begin{proposition} \label{prop:square sum of BesselJ}
		Let $\nu \in \intervalco{-1/2}{\infty}$ and $r \in \intervaloo{0}{\infty}$. 
		Then we have
		\begin{equation} \label{eq:square sum of BesselJ}
			\frac{d}{dr} \mleft( r ( \BesselJ{\nu}(r)^2 + \BesselJ{\nu+1}(r)^2 ) \mright) = (2\nu+1) ( \BesselJ{\nu}(r)^2 - \BesselJ{\nu+1}(r)^2 ) .
		\end{equation}
		As a consequence, we have
		\begin{equation} \label{eq:square sum of BesselJ Tnu}
			\frac{d}{ds} ( \Tnu{\nu} f(s) + \Tnu{\nu+1} f(s) ) = \frac{2\nu+1}{s} ( \Tnu{\nu} f(s) - \Tnu{\nu+1} f(s) )
		\end{equation}
		for every $f \in \Lpnunu{1}{2\nu+1}{0}$.
	\end{proposition}
	\begin{proof}[Proof of Proposition \ref{prop:square sum of BesselJ}]
		It is well known that 
		\begin{align}
			\frac{d}{dr} \BesselJ{\nu}(r) &= - \BesselJ{\nu+1}(r) + \frac{\nu}{r} \BesselJ{\nu}(r) , \\
			\frac{d}{dr} \BesselJ{\nu+1}(r) &= \BesselJ{\nu}(r) - \frac{\nu+1}{r} \BesselJ{\nu + 1}(r)
		\end{align}
		hold (see \cite[\href{https://dlmf.nist.gov/10.6.E2}{10.6.2}]{DLMF}, \cite[8.472]{GR2014}).
		Using these, we get
		\begin{align}
			\frac{d}{dr} \mleft( r \BesselJ{\nu}(r)^2 \mright)
			&= \BesselJ{\nu}(r)^2 + 2 r \BesselJ{\nu}(r) \mleft( - \BesselJ{\nu+1}(r) + \frac{\nu}{r} \BesselJ{\nu}(r) \mright) \\
			&= (2\nu+1) \BesselJ{\nu}(r)^2 - 2 r \BesselJ{\nu}(r) \BesselJ{\nu+1}(r)
			\shortintertext{and}
			\frac{d}{dr} \mleft( r \BesselJ{\nu+1}(r)^2 \mright)
			&= \BesselJ{\nu+1}(r)^2 + 2 r \BesselJ{\nu+1}(r) \mleft( \BesselJ{\nu}(r) - \frac{\nu+1}{r} \BesselJ{\nu + 1}(r) \mright) \\
			&= - (2\nu+1) \BesselJ{\nu+1}(r)^2 + 2 r \BesselJ{\nu}(r) \BesselJ{\nu+1}(r) ,
		\end{align}
		so that \eqref{eq:square sum of BesselJ} holds. 
		\eqref{eq:square sum of BesselJ Tnu} is immediate by differentiating under the integral sign (which is justified by the estimate \eqref{eq:BesselJ upper bound}).
	\end{proof}
	We remark that $\Tnu{\nu} f$ and $\Tnu{\nu+1} f$ can be non-differentiable in Proposition \ref{prop:square sum of BesselJ}, even though $\Tnu{\nu} f + \Tnu{\nu+1} f$ is differentiable. For example, let $\nu = 1/2$ and $f(r) = (1 - \cos{r})/r^2$.
	In this case, one can show that
	\begin{align}
		\Hankelnu{1/2}f(s) 
		&= \begin{cases}
			(\pi / 2)^{1/2} / s , & s \in \intervaloo{0}{1} , \\
			(\pi / 2)^{1/2} / 2 , & s = 1 , \\
			\mathrlap{0,}\hphantom{\pi/2 - \pi/(12s^2),} & s \in \intervaloo{1}{\infty} , 
		\end{cases}\\
		\Tnu{1/2}f(s) 
		&= \begin{cases}
			\pi s , & s \in \intervaloc{0}{1/2} , \\
			\mathrlap{\pi / 2,}\hphantom{\pi/2 - \pi/(12s^2),} & s \in \intervalco{1/2}{\infty} , 
		\end{cases}\\
		\Tnu{3/2}f(s) 
		&= \begin{cases}
			\pi s / 3, & s \in \intervaloc{0}{1/2} , \\
			\pi/2 - \pi/(12s^2) , & s \in \intervalco{1/2}{\infty} .
		\end{cases}
	\end{align}
	See \citet[3.741.2, 6.672.5]{GR2014} for these identities.
	Clearly, $\Tnu{1/2}f$ and $\Tnu{3/2}f$ are not differentiable at $s = 1/2$.
	Nevertheless, 
	\begin{equation}
		\Tnu{1/2}f(s) + \Tnu{3/2}f(s) = \begin{cases}
			4 \pi s / 3, & s \in \intervaloc{0}{1/2} , \\
			\pi - \pi/(12s^2) , & s \in \intervalco{1/2}{\infty} 
		\end{cases}
	\end{equation}
	is differentiable on $\intervaloo{0}{\infty}$.

\end{document}